\documentclass[hidelinks]{article}

\usepackage[utf8]{inputenc}
\usepackage{geometry}
\usepackage{amsmath}
\usepackage{amssymb}
\usepackage{amsfonts}
\usepackage{amsthm}
\usepackage{hyperref}
\usepackage{cleveref}
\usepackage{latexsym}
\usepackage{array}
\usepackage{graphicx}
\usepackage{caption}
\usepackage{algorithmic}
\usepackage{float}
\usepackage{algorithm}
\usepackage{verbatim}
\usepackage{xcolor}

\usepackage{graphicx,epstopdf}
\usepackage[caption=false]{subfig}

\usepackage{authblk}

\newtheorem{theorem}{Theorem}
\newtheorem{definition}{Definition}
\newtheorem{lemma}{Lemma}
\newtheorem{remark}{Remark}

\title{High order asymptotic preserving scheme for linear kinetic equations with diffusive scaling}
\date{}

\author{Megala Anandan\thanks{\noindent Indian Institute of Science, C.V. Raman Road, 560012, Bangalore, India
  (megalaa@iisc.ac.in).}, 
Benjamin Boutin{\footnote{
Univ Rennes, CNRS, IRMAR UMR 6625, 35000 Rennes, France 
  (benjamin.boutin@univ-rennes.fr).}}, Nicolas Crouseilles{\footnote{
Univ Rennes, CNRS, IRMAR UMR 6625 \& centre Inria de l’Université de Rennes (MINGuS) \& ENS Rennes, France (nicolas.crouseilles@inria.fr).}  }}

\begin{document}

\maketitle

\begin{abstract}
    In this work, high order asymptotic preserving schemes are constructed and analysed for kinetic equations under a diffusive scaling. The framework enables to consider different cases: the diffusion equation, the advection-diffusion equation and the presence of inflow boundary conditions. Starting from the micro-macro reformulation of the original kinetic equation, high order time integrators are introduced. This class of numerical schemes  enjoys the Asymptotic Preserving (AP) property for arbitrary initial data and degenerates when $\epsilon$ goes to zero into a high order scheme which is implicit for the diffusion term, which makes it free from the usual diffusion stability condition. The space discretization is also discussed and high order methods are also proposed based on classical finite differences schemes. The Asymptotic Preserving property is analysed and numerical results are presented to illustrate the 
    properties of the proposed schemes in different regimes. 
\end{abstract}

\section{Introduction}
In this work, we are concerned with the numerical approximation of linear kinetic transport equations in a diffusive scaling. Such models are 
widely used in applications such as rarefied gas dynamics, neutron transport, and radiative transfer. 
Due to the presence of a small parameter $\epsilon$ 
(which is the normalized mean free path of the particles), 
standard schemes suffer from a severe restriction on 
the numerical parameters, making the simulations 
very costly. In the last decades, the so-called Asymptotic-Preserving (AP) schemes have been proposed to make possible the numerical passage between the micro and macro scale \cite{ap_jin, ap_jin2}. Indeed, these AP schemes are uniformly stable and degenerate when $\epsilon\to 0$ to a scheme which is consistent with the asymptotic diffusion model. This makes them very attractive to deal with multi-scale phenomena as an alternative to domain decomposition approaches. 

The goal of this work is to design high order in time AP 
schemes for collisional kinetic equations in the diffusive scaling. Several works can be found in the literature 
on this topic \cite{jin_pareschi_toscani, ap_jin, ap_jin2, klar1, klar2, klar3, klar3, pareschi, dimarco, lafitte, larsen, doi:10.1137/07069479X, Lemou2010RelaxedMS, PENG2020109485}. Our work is based on a micro-macro decomposition as introduced in \cite{doi:10.1137/07069479X} where the unknown $f$ of the stiff kinetic equation is split into an equilibrium part $\rho$ plus a remainder $g$. A micro-macro model (equivalent to the original kinetic one) satisfied by $\rho$ ad $g$ can be derived. This micro-macro strategy turns out to be the starting point of several numerical approximations in phase space (using particles method, Discontinuous Galerkin method or low rank approximation \cite{crous_lemou, crestetto, JANG2015199, low_rank1, low_rank2}). In addition, a suitable first order semi-implicit time discretization of the micro-macro model is used as in \cite{doi:10.1137/07069479X} for which however the asymptotic diffusion equation is solved explicitly. This drawback is overcome 
following \cite{Lemou2010RelaxedMS, crous_lemou, crestetto} in which the AP scheme degenerates into an implicit treatment of the diffusion equation. This improvement 
enables to get a numerical scheme which is asymptotically free from the usual parabolic condition. 

The derivation of high order in time AP schemes 
for stiff kinetic problem has been performed by several authors \cite{dimarco_exprk, dimarco_imexrk, dimarco_multistep, albi_multistep, JANG2015199}
using the so-called high order IMEX methods \cite{imex1, imex2, pareschi_russo}. 
In this work, a family of high order IMEX schemes is proposed for linear  collisional kinetic equations in the diffusive scaling, which degenerates when $\epsilon \to 0$ 
to a high order IMEX scheme for the diffusion equation. From the first order semi-implicit AP numerical scheme \cite{crous_lemou}, 
the family of high order schemes proposed in this work is  obtained using globally stiffly accurate high order IMEX 
Runge-Kutta methods, namely type A and type CK \cite{dimarco_imexrk, JANG2015199}. 

In addition to the standard diffusion scaling, 
we also consider two other examples that enter in our framework. First, we consider a modification of the collision operator that enables to derive a transport-diffusion asymptotic model  \cite{jin_pareschi_toscani, JANG2015199}. Second,  we discuss half moments micro-macro decomposition which naturally incorporates the incoming boundary conditions \cite{doi:10.1137/120865513}. 

Lastly, we address the space discretization 
in order to get a fully high order solver of the kinetic equation. 
High order space approximation 
based on finite difference methods is considered. Staggered or non-staggered strategies are proposed to achieve high order accuracy in space. 

The paper is organized as follows. First in Section~\ref{sec:2}, the kinetic and asymptotic diffusion models are introduced. 
Then in Section~\ref{sec:3}, high order time integrators  (using globally stiffly accurate IMEX Runge-Kutta temporal discretization) 
are proposed, and their AP property in the diffusive limit is addressed in Section~\ref{sec:4}. Section~\ref{sec:5} is devoted to the space approximation. 
In Section~\ref{sec:6}, we discuss some extensions to other collision operators and to half moments. In Section~\ref{sec:7}, numerical results are presented, illustrating high order accuracy and 
the main properties of the schemes.

\section{Kinetic equation, diffusion limit and micro-macro decomposition}
\label{sec:2}
In this section, we introduce the kinetic model in the diffusive scaling, and recall the asymptotic limit. 
Then, the micro-macro decomposition is performed to derive the micro-macro model which serves as a basis for the 
numerical developments. 

\subsection{Linear kinetic equation with diffusive scaling}
\label{Subsec:Intro kinetic eqn}
Let $\Omega \subset \mathbb{R}^d$ be the position space and $V \subseteq \mathbb{R}^d$ be the velocity space with measure $d\mu(v)$. We consider the linear kinetic equation with diffusive scaling,
\begin{equation}
\label{BGK eqn}
    \partial_t f + \frac{1}{\epsilon} v \cdot \nabla_x f = \frac{1}{\epsilon^2}Lf, \quad (t,x,v) \in \mathbb{R}^+ \times \Omega \times V
\end{equation}
where $f(t,x,v) \in \mathbb{R}$ is the distribution function  (depending on time $t\in\mathbb{R}^+$,   space $x\in\Omega\subset \mathbb{R}^d$  and  velocity $v\in V \subset \mathbb{R}^d$) and $\epsilon>0$ measures the dimensionless mean free path of particles or the inverse of  relaxation time. We consider the initial condition,
\begin{equation}
    f(0,x,v)=f^{\textsf{init}}(x,v), \quad (x,v) \in \Omega \times V
\end{equation}
and boundary conditions are imposed in space. 
In this work, we will consider periodic boundary conditions 
or inflow boundary conditions. 
The linear collision operator $L$ in \cref{BGK eqn} acts only on the velocity dependence of $f$, and it relaxes the particles to an equilibrium $M(v)$ which is positive and even. We denote for all velocity dependent distribution functions $h$, 
\begin{equation}
\label{def_crochet}
    \left<h\right>_V =\frac{\int_{V} h(v) \  d\mu }{\int_{V} M(v) \ d\mu}.
\end{equation}
In particular, we obtain $\left<M\right>_V =1$ and $\left<v M\right>_V =0$. Further, the  operator $L$ is non-positive and self-adjoint in $L^2\left( V, M^{-1} d\mu \right)$, with the following null space and range:
\begin{equation}
\label{L Nullspace}
    \mathcal{N}(L)=\{f:f\in \mathsf{Span}\left( M \right) \}, 
 \;\;    \mathcal{R}(L)= \left( \mathcal{N}(L) \right)^\perp = \{ f: \left<f\right>_V = 0 \} .
\end{equation}
Therefore, $L$ is invertible on $\mathcal{R}(L)$ and we denote its pseudo-inverse by $L^{-1}$. We also assume that $L$ is invariant under orthogonal transformations of $\mathbb{R}^d$. 

\subsection{Diffusion limit}
In the limit $\epsilon \rightarrow 0$, it is seen from \cref{BGK eqn} that $f\to f_0$ where $f_0$ belongs to ${\cal N}(L)$. Thus, $f_0=\rho (t,x) M$ where $f_0$ solves $Lf_0=0$ and where the limiting density $\rho$ is the solution 
of the asymptotic diffusion equation. To derive the diffusion equation, a Chapman-Enskog expansion has to be performed to get $f=f_0+\epsilon L^{-1}(vM) \cdot \nabla_x \rho + {\cal O}(\epsilon^2)$.
Integrating with respect to the velocity variable enables to get the diffusion limit 
\begin{equation}
\label{Diff eqn}
    \partial_t \rho - \nabla_x \cdot \left( \kappa \nabla_x \rho \right) = 0 \text{ with } \kappa= - \left< v \otimes L^{-1}(vM) \right>_V > 0. 
\end{equation}

\subsection{Micro-macro decomposition}
In this part, we derive a micro-macro model which is equivalent to \eqref{BGK eqn}, and this is the model that 
will be discretized in the next sections. 
First, we consider the standard micro-macro decomposition 
of the unknown $f$  \cite{doi:10.1137/07069479X,Lemou2010RelaxedMS},
\begin{equation}
\label{mic-mac}
    f=\rho M + g, \;\; \mbox{ with } \;\; 
   \rho(t,x)=\left< f \right>_V \text{ and } \left< g \right>_V = 0.
\end{equation}
We introduce the orthogonal projector $\Pi$ in $L^2\left( V, M^{-1} d\mu \right)$ onto $\mathcal{N}(L)$: $\Pi h=\left<h\right>_V M$, 
which will be useful to derive the micro-macro model. 
Substituting \cref{mic-mac} into \cref{BGK eqn} and applying successively $\Pi$ and $(I-\Pi)$ enables to get 
the micro-macro model satisfied by $(\rho, g)$
\begin{gather}
\label{macro std}
    \partial_t \rho + \frac{1}{\epsilon} \nabla_x \cdot \left< v   g\right>_V = 0,  \\
\label{micro std}
    \partial_t g + \frac{1}{\epsilon} \left( I-\Pi \right) \left( v\cdot\nabla_x g \right) + \frac{1}{\epsilon}vM\cdot\nabla_x\rho = \frac{1}{\epsilon^2}Lg.
\end{gather}
Initial conditions for macro and micro equations become  
\begin{gather}
\label{ID rho}
    \rho(0,x)=\rho^{\textsf{init}}(x)=\left< f^{\textsf{init}}(x, \cdot) \right>_V,  \\
\label{ID g1}
    g(0,x,v)=g^{\textsf{init}}(x,v)= f^{\textsf{init}}(x, v) - \rho^{\textsf{init}}(x) M(v), 
\end{gather}
whereas the boundary conditions for $\rho$ and $g$ become periodic if $f$ is periodic. 
From the micro part \eqref{micro std}, 
a Chapman-Enskog expansion of $g$ can be performed to get 
$$
    g=-\epsilon\left( \epsilon^2 \partial_t - L \right)^{-1}\Big(  (I-\Pi)\left( v \cdot \nabla_x g \right) + vM\cdot \nabla_x\rho \Big) = \epsilon L^{-1}(vM) \cdot \nabla_x\rho + {\cal O}(\epsilon^2), 
$$
under some suitable smoothness assumptions. 
Inserting this expression in \cref{macro std} 
leads to \cref{Diff eqn} in the limit $\epsilon \rightarrow 0$. 

\section{Time integrators}
\label{sec:3}
In this part, we present the family  of high order time integrators for the micro-macro model \eqref{macro std}-\eqref{micro std}. We will keep the phase space variables continuous to ease the reading. 
We first recall the first order temporal scheme which leads to the implicit treatment of the asymptotic diffusion model before introducing the high order version.

\subsection{First order accurate time integrator}
Given $\rho^n, g^n$ that approximate $\rho, g$ at time $t=n\Delta t$, we obtain the solution $\rho^{n+1}, g^{n+1}$ from the following time integration of \cref{macro std,micro std} respectively. We use the following 
first order implicit-explicit (IMEX) strategy to attain the asymptotic preserving property
\begin{gather}
\label{macro std 1st}
    \rho^{n+1}=\rho^n - \frac{\Delta t}{\epsilon} \nabla_x \cdot \left< v g^{n+1}\right>_V, \\
\label{micro std 1st}
    g^{n+1} = g^n - \frac{\Delta t}{\epsilon} \left( I-\Pi \right) \left( v\cdot\nabla_x g^n \right) - \frac{\Delta t}{\epsilon}  vM\cdot\nabla_x\rho^{n+1} + \frac{\Delta t}{\epsilon^2} Lg^{n+1}.
\end{gather}
 Let us observe that this scheme is different from the IMEX strategies employed in \cite{doi:10.1137/07069479X,JANG2015199}, due to our implicit treatment of density gradient in micro equation 
 \eqref{micro std 1st} and fully implicit treatment of the macro equation. This strategy enables us to get an implicit scheme for diffusion equation in the limit $\epsilon \rightarrow 0$. \\
Although the macro equation is treated in a fully implicit manner, $\rho^{n+1}$ and $g^{n+1}$ can be updated using \cref{macro std 1st,micro std 1st} in an explicit manner. From \cref{micro std 1st}, we get
\begin{equation}
\label{g_imp_first}
    g^{n+1}= \left( \epsilon^2 I - \Delta t L \right)^{-1} \left( \epsilon^2 g^n - \epsilon \Delta t \left( I-\Pi \right) \left( v\cdot\nabla_x g^n \right) - \epsilon \Delta t  vM\cdot\nabla_x\rho^{n+1} \right). 
\end{equation}
Inserting this in \cref{macro std 1st}, we obtain 
the following implicit scheme for the macro unknown 
\begin{equation*}
\rho^{n+1}\! =\! \rho^n -\Delta t \nabla_x \cdot \langle v \left( \epsilon^2 I - \Delta t L \right)^{-1} \!\!\left( \epsilon g^n \!-\! \Delta t \left( I-\Pi \right) \left( v\cdot\nabla_x g^n \right) \!-\!\Delta t vM \!\cdot\! \nabla_x \rho^{n+1}  \right)\rangle_{\tiny V}\!, 
\end{equation*}
or, expressing $\rho^{n+1}$ as quantities at iteration $n$  
\begin{equation*}
\label{rho_imp_first}
\rho^{n+1} = \left( I -\Delta t^2 \nabla_x \cdot \left(\mathcal{D}_{\epsilon, \Delta t} \nabla_x \right) \right)^{-1} 
\biggl( \rho^n -\Delta t \nabla_x \!\cdot\! \biggl< v \left( \epsilon^2 I - \Delta t L \right)^{-1}  \left( \epsilon g^n - \Delta t \left( I-\Pi \right) ( v\!\cdot\!\nabla_x g^n )  \right)\biggr>_V  \biggr) 
\end{equation*}
with $\mathcal{D}_{\epsilon, \Delta t} = \langle v \otimes \left( \epsilon^2 I -  \Delta t L \right)^{-1} \left(vM\right) \rangle_V$. 
Thanks to this time integrator, $\rho^{n+1}$ can be updated by inverting a diffusion type operator. Following this, $g^{n+1}$ can be found explicitly from the knowledge of $\rho^{n+1}$. 
This first order scheme introduced in \cite{Lemou2010RelaxedMS, crous_lemou} 
is the basis of the high order scheme 
presented below. 

\subsection{High order accurate time integrators}
Following previous works \cite{dimarco_imexrk, JANG2015199,doi:10.1137/110842855}, we will consider globally stiffly accurate (GSA) IMEX Runge-Kutta (RK) schemes to construct high order  time integrators for the micro-macro model \cref{macro std,micro std}. An IMEX RK scheme is represented using the double Butcher tableau \cite{imex1, imex2}
\begin{equation}
\label{gen Butcher tableau}
\begin{tabular}{p{0.25cm}|p{0.5cm}}
\centering $\Tilde{c}$ & \centering $\Tilde{A}$ \cr 
\hline
 & \centering \vspace{-0.25em}$\Tilde{b}^T$
\end{tabular} \quad \quad
\begin{tabular}{p{0.25cm}|p{0.5cm}}
\centering $c$ & \centering $A$ \cr 
\hline
 &  \vspace{-0.25em}\centering $b^T$ 
\end{tabular}
\end{equation}
where $\Tilde{A}=(\Tilde{a}_{ij})$ and $A=(a_{ij})$ 
are $s \!\times \!s$ matrices which correspond to the explicit and implicit parts of the scheme ($A$ and $\Tilde{A}$ respectively are lower triangular and strictly lower triangular matrices). The coefficients $\Tilde{c}$ and $c$ are given by  $\Tilde{c}_i=\sum_{j=1}^{i-1} \Tilde{a}_{ij}$, $c_i=\sum_{j=1}^{i} a_{ij}$, and the vectors $\Tilde{b}=(\Tilde{b}_j)$ and $b=(b_j)$ give quadrature weights that combine the stages. For GSA IMEX RK scheme, we have 
\begin{equation}
\label{GSA}
    c_s=\Tilde{c}_s=1 \text{ and } a_{sj}=b_j, \Tilde{a}_{sj}=\Tilde{b}_j, \;\; \forall j \in \{1,2..,s\}.  
\end{equation}
An IMEX RK method is type A if the matrix $A$ is invertible, and it is type CK if the first row of matrix $A$ has zero entries and the square sub-matrix formed by excluding the first column and row of $A$ is invertible. In the special case where the first column of $A$ has zero entries, the scheme is said to be of type CK-ARS. The reader is referred to \cite{dimarco_imexrk} for more details. In this work, we employ both type A and CK-ARS schemes. \\  
The first order GSA IMEX RK scheme employed in \cref{macro std 1st,micro std 1st} follows the type CK-ARS double Butcher tableau (known as ARS$(1,1,1)$),
\begin{equation}
\begin{tabular}{p{0.25cm}|p{0.25cm}p{0.25cm}}
\centering $0$ & $0$ & $0$ \cr 
\centering $1$ & $1$ & $0$ \cr 
\hline
 & $1$ & $0$
\end{tabular} \quad \quad
\begin{tabular}{p{0.25cm}|p{0.25cm}p{0.25cm}}
\centering $0$ & $0$ & $0$ \cr 
\centering $1$ & $0$ & $1$ \cr 
\hline
 & $0$ & $1$
\end{tabular} 
\end{equation}
We now use the general IMEX RK scheme from \eqref{gen Butcher tableau} with GSA property \cref{GSA} for obtaining high order accurate time integration of macro and micro \cref{macro std,micro std} respectively. We introduce the following notations in the presentation of our scheme. 
\begin{gather}
    \label{def_not_T}
    \mathcal{T}h^{(k)}=\left( I-\Pi \right) \left( v\cdot\nabla_x h^{(k)} \right), \\ 
    \label{def_not_D}
    \mathcal{D}^{(j)}_{\epsilon, \Delta t} = \left< v \otimes \left( \epsilon^2 I - a_{jj} \Delta t L \right)^{-1} \left(vM\right) \right>_V, \\ 
    \label{def_not_I}
    \mathcal{I}^{(j)}_{\epsilon, \Delta t} = \left( \epsilon^2 I - a_{jj} \Delta t L \right)^{-1}. 
\end{gather}
 We will construct high order IMEX RK schemes following the first order guidelines (fully implicit treatment of macro equation, implicit treatment of density gradient and relaxation terms and explicit treatment of transport term in micro equation).
Given $\rho^n, g^n$ that approximate $\rho, g$ at time $t=n\Delta t$, we obtain the internal RK stage values $\rho^{(j)}$ and $g^{(j)}$, $j=1,\dots,s$ as
\begin{gather}
\label{macro std high stages}
    \rho^{(j)}=\rho^n - \sum_{k=1}^j a_{jk} \frac{\Delta t}{\epsilon} \nabla_x \cdot  \left< v g^{(k)}\right>_V, \\
\label{micro std high stages}
    g^{(j)} = g^n - \sum_{k=1}^{j-1} \Tilde{a}_{jk} \frac{\Delta t}{\epsilon} \mathcal{T}g^{(k)} - \sum_{k=1}^j a_{jk} \frac{\Delta t}{\epsilon}  vM\cdot\nabla_x\rho^{(k)} + \sum_{k=1}^j a_{jk} \frac{\Delta t}{\epsilon^2} Lg^{(k)}, 
\end{gather}
where, as usual, the summation $\sum_{k=1}^{j-1}$ in the explicit term  is zero for $j=1$. \\
Although the expressions above are implicit, the stage values $\rho^{(1)}$, $g^{(1)}$ can be found in an explicit manner by using the known quantities $\rho^n, g^n$, and the stage values $\rho^{(j)}$, $g^{(j)}, \ \forall j \in \{2,3,\dots,s \}$ can be found explicitly from $\rho^n, g^n$ and the previous stage values $\rho^{(l)},g^{(l)}, \ \forall l \in \{ 1,2,\dots,j-1\}$. 
Indeed, proceeding similarly as for the first order scheme, we get the following expression of $g^{(j)}, j=1, \dots, s$ from \cref{micro std high stages}, 
\begin{equation}
\label{Exp for g_std high}
    g^{(j)}= \mathcal{I}^{(j)}_{\epsilon, \Delta t} \left( \epsilon^2 g^n - \epsilon \sum_{k=1}^{j-1} \Tilde{a}_{jk} \Delta t \mathcal{T} g^{(k)}  - \epsilon \sum_{k=1}^{j} a_{jk} \Delta t  vM\cdot\nabla_x\rho^{(k)} + \sum_{k=1}^{j-1} a_{jk} \Delta t Lg^{(k)} \right). 
\end{equation}
Further, we write \cref{macro std high stages} by splitting the summation on $k$ as 
\begin{equation*}
    \rho^{(j)}=\rho^n - \sum_{k=1}^{j-1} a_{jk} \frac{\Delta t}{\epsilon} \nabla_x \cdot \left< v g^{(k)}\right>_V - a_{jj} \frac{\Delta t}{\epsilon} \nabla_x \cdot \left< v g^{(j)}\right>_V, 
\end{equation*}
and inserting \cref{Exp for g_std high} in the last term leads to the update of  $\rho^{(j)}$ for $j=1,\dots,s$
\begin{eqnarray}
\label{Exp for rho_std high}
    \rho^{(j)} &=& \left( I -a_{jj}^2 \Delta t^2  \nabla_x \cdot \left( \mathcal{D}^{(j)}_{\epsilon, \Delta t} \nabla_x \right) \right)^{-1}  
    \Bigg( \rho^n - \sum_{k=1}^{j-1} a_{jk}\frac{\Delta t}{\epsilon} \nabla_x \cdot \left< v g^{(k)}\right>_V  \\
    &&- a_{jj} \Delta t \nabla_x \cdot \Big< v \mathcal{I}^{(j)}_{\epsilon, \Delta t}
    \bigl( \epsilon g^n - \sum_{k=1}^{j-1} \Tilde{a}_{jk} \Delta t \mathcal{T} g^{(k)} - \sum_{k=1}^{j-1} a_{jk} \Delta t  vM\cdot\nabla_x\rho^{(k)}  + \frac{1}{\epsilon} \sum_{k=1}^{j-1} a_{jk} \Delta t Lg^{(k)}\bigr)\Big>_V  \Bigg), \nonumber 
\end{eqnarray}
where the definition of $\mathcal{T},  \mathcal{D}^{(j)}_{\epsilon, \Delta t}$ and 
$ \mathcal{I}^{(j)}_{\epsilon, \Delta t}$ 
are given by \cref{def_not_T,def_not_D,def_not_I}.  
After this reformulation, $\rho^{(j)}$ can be computed from \eqref{Exp for rho_std high} by inverting a linear elliptic type problem and following this, $g^{(j)}$ can be found from \cref{Exp for g_std high}.  
The GSA property in \cref{GSA} guarantees that the solution at $t^{n+1}=(n+1)\Delta t$ is same as the last RK stage values, that is, $\rho^{n+1}=\rho^{(s)}$ and $g^{n+1}=g^{(s)}$. 

\section{Asymptotic preserving property}
\label{sec:4}
In this section, we show that the time integrated scheme \eqref{Exp for rho_std high}-\eqref{Exp for g_std high} becomes a consistent scheme for the diffusion equation \eqref{Diff eqn} in the limit $\epsilon \rightarrow 0$. We will discuss the asymptotic preserving property for both CK-ARS and type A time integrators as performed in \cite{dimarco_imexrk} for the fluid limit.  
First, we recall the definition of well-prepared initial data in our context. 
\begin{definition}[Well-prepared initial data]
\label{Def WP ID}
$\!\!\!\!\!$  The initial data $\rho(0,x)$ and $g(0,x,v)$ in \cref{ID rho,ID g1} are said to be well-prepared if $g(0,x,v)=O(\epsilon)$.
\end{definition}
\begin{lemma}
\label{Lem AP g}
    Assume that $\epsilon$ is sufficiently small.
    Let $\Tilde{a}_{jk}$ and $a_{jk}$ be the coefficients of the RK method \eqref{gen Butcher tableau} applied to the scheme \eqref{macro std high stages}-\eqref{micro std high stages}. 
    Then, the following holds:
    \begin{enumerate}
        \item CK-ARS case: 
        If $g^n=O(\epsilon)$, then  $g^{(1)}=g^n=O(\epsilon)$ and \\
        $g^{(j)}=\epsilon L^{-1} (vM)\cdot\nabla_x \rho^{(j)}+O\left(\epsilon^2\right), \ \forall j \in \{2,\dots,s \}$. 
        \item Type A case: 
        $g^{(j)}=\epsilon L^{-1} (vM)\cdot\nabla_x \rho^{(j)}+O\left(\epsilon^2\right), \ \forall j \in \{1,\dots,s \}$.
    \end{enumerate}
\end{lemma}

\begin{proof}
    Let $j\in\{1,\ldots,s\}$ such that $a_{jj}\neq 0$. Observe that the operator $\mathcal{I}^{(j)}_{\epsilon, \Delta t}$ defined in~\eqref{def_not_I} admits, for small $\epsilon$, the following expansion:
    \begin{equation}
    \label{exp_Ij}
        \mathcal{I}^{(j)}_{\epsilon, \Delta t} 
        = -(a_{jj}\Delta t L)^{-1} + O(\epsilon^2).
    \end{equation}

    Consider now an A-type time integrator, so with $a_{jj}\neq 0$ for any $j\in\{1,\ldots,s\}$, and assume $g^n = O(1)$. From~\eqref{Exp for g_std high} and the previous expansion, we obtain
    \begin{align*}
        g^{(1)}
        = -(a_{11}\Delta t L)^{-1} \left[ -\epsilon a_{11}\Delta t vM\cdot \nabla_x \rho^{(1)}\right] + O(\epsilon^2)
        = \epsilon L^{-1}(vM)\cdot\nabla_x \rho^{(1)} + O(\epsilon^2).
    \end{align*}
    Now, the proof is performed by induction on $j\in\{2,\ldots,s\}$ assuming that for any $k\in\{1,\ldots,j-1\}$, $g^{(k)} = \epsilon L^{-1}(vM)\cdot\nabla_x \rho^{(k)} + O(\epsilon^2)$. In particular $g^{(k)}=O(\epsilon)$ and the formula~\eqref{Exp for g_std high} has therefore the following expansion:
    \begin{equation*}
        g^{(j)}
        = -(a_{jj}\Delta t L)^{-1} \left[ O(\epsilon^2)  - \epsilon \sum_{k=1}^{j} a_{jk} \Delta t  vM\cdot\nabla_x\rho^{(k)} + \sum_{k=1}^{j-1} a_{jk} \Delta t Lg^{(k)}  \right] +  O(\epsilon^2).
    \end{equation*}
    Inserting the induction hypothesis in the last sum, most of the terms in the two sums eliminate so that finally $g^{(j)} = \epsilon L^{-1}(vM)\cdot\nabla_x \rho^{(j)} + O(\epsilon^2)$.

    The case of a CK-ARS time integrator is slightly different. First $a_{11}=0$ so that $g^{(1)}=g^n = O(\epsilon)$ by the particular well-prepared assumption. Now $a_{22}\neq 0$ and~\eqref{Exp for g_std high} has the following expansion for $j=2$:
    \begin{align*}
        g^{(2)}\! =\! -(a_{22}\Delta t L)^{-1} \!\left[ O(\epsilon^2) \!-\! \epsilon a_{22} \Delta t vM\cdot\nabla_x\rho^{(2)} \!\right] \!+\! O(\epsilon^2)
        \!=\! \epsilon L^{-1}(vM)\!\cdot\!\nabla_x \rho^{(2)} + O(\epsilon^2).
    \end{align*}
    Again, the proof is by induction on $j\in\{3,\ldots,s\}$ assuming for any $k\in\{2,\ldots,j-1\}$, $g^{(k)} = \epsilon L^{-1}(vM)\cdot\nabla_x \rho^{(k)} + O(\epsilon^2)$. The same computation as above is available since $g^{(1)}=O(\epsilon)$. One has (note that $a_{j1}=0$ for any $j$ so that the sums start at $k=2$):
    \begin{align*}
        g^{(j)}
        &= -(a_{jj}\Delta t L)^{-1} \left[ O(\epsilon^2)  - \epsilon \sum_{k=2}^{j} a_{jk} \Delta t  vM\cdot\nabla_x\rho^{(k)} + \sum_{k=2}^{j-1} a_{jk} \Delta t Lg^{(k)}  \right] +  O(\epsilon^2)\\
        &= \epsilon L^{-1}(vM)\cdot\nabla_x \rho^{(j)} + O(\epsilon^2).
    \end{align*}
\end{proof}
Due to the GSA property of both time integrators, we have $g^{n+1}=g^{(s)}=\epsilon L^{-1} (vM)\cdot\nabla_x \rho^{(s)}+O\left(\epsilon^2\right)=\epsilon L^{-1} (vM) \cdot\nabla_x \rho^{n+1}+O\left(\epsilon^2\right)$ for sufficiently small $\epsilon$. Thus, the following are evident from \cref{Lem AP g}: 
\begin{enumerate}
    \item For type CK-ARS, if the initial data is well-prepared (that is, $g^0=O(\epsilon)$), then $g^n=O(\epsilon), \ \forall n>0$.
    \item For type A, if the initial data is such that $g^0=O(1)$, then $g^n=O(\epsilon), \ \forall n>0$.
\end{enumerate}
As observed in \cite{dimarco_imexrk}, the initial data does not need to be well-prepared for type A time integrator, unlike type CK-ARS, to ensure AP property. 
\begin{theorem} 
$\!\!$Consider the scheme \eqref{macro std high stages}-\eqref{micro std high stages} 
approximating the macro-micro model \eqref{macro std}-\eqref{micro std},
with the RK method \eqref{gen Butcher tableau} of type A or of type CK-ARS (with well-prepared initial data $g^0=O(\epsilon)$).  
Then in the limit $\epsilon \rightarrow 0$, the scheme \eqref{macro std high stages}-\eqref{micro std high stages} 
degenerates to the following scheme for the diffusion equation
\begin{equation}
\label{rho_j thm AP}
        \rho^{(j)}=\rho^n + \sum_{k=1}^j a_{jk} \Delta t \nabla_x \cdot \left(\kappa \nabla_x \rho^{(k)}\right), \ \forall j=1, \dots, s, \;   \kappa\!=\! - \left< v \otimes L^{-1}\!(vM) \right>_V. 
    \end{equation}
\end{theorem}
\begin{proof}
    Corresponding to each case (CK-ARS or type A), we have the following:
    \begin{description}
        \item[Type CK-ARS] Assumptions in criterion 1 of Lemma \ref{Lem AP g}  are satisfied, and its implications can be utilised. Hence, inserting $g^{(\ell)}=\epsilon L^{-1} (vM) \cdot\nabla_x \rho^{(\ell)}+O(\epsilon^2), \ \forall \ell \in \{2,3,..,s \}$ into \cref{macro std high stages}, we get (recall that $a_{j1}=0$) 
        \begin{align*}
            \rho^{(j)}
            &=\rho^n - \frac{\Delta t}{\epsilon} \sum_{k=2}^j a_{jk} \nabla_x \cdot \left< v \epsilon L^{-1} (vM) \cdot \nabla_x \rho^{(k)} \right>_V + O(\epsilon), \\
            &=\rho^n - \Delta t \sum_{k=2}^j a_{jk} \nabla_x \cdot \left( \left< v \otimes L^{-1} (vM) \right>_V \nabla_x \rho^{(k)} \right) + O(\epsilon).
        \end{align*}
        
        \item[Type A] Assumptions in criterion 2 of Lemma \ref{Lem AP g} are satisfied, and its implications can be utilised. Hence, inserting $g^{(\ell)}=\epsilon L^{-1} (vM) \cdot\nabla_x \rho^{(\ell)}+O(\epsilon^2), \ \forall \ell \in \{1,2,..,s \}$ into \cref{macro std high stages}, we get the required result by following the same simplification as before. The only difference is that here $\sum_{k=1}^j$ instead of $\sum_{k=2}^j$.
    \end{description}
\end{proof}

\begin{remark}
\label{rem_nwp}
     For type CK-ARS, if the initial data is not well-prepared, computing $g^{(2)}$ from \eqref{micro std high stages} involves $\epsilon \frac{\Tilde{a}_{21}}{a_{22}} L^{-1} (I-\Pi)(v\cdot\nabla_x g^{(1)})$ which is not of $O(\epsilon^2)$.  Thus,
     \begin{equation*}
         g^{(2)}=\epsilon \frac{\Tilde{a}_{21}}{a_{22}} L^{-1} (I-\Pi)(v\cdot\nabla_x g^{(1)}) + \epsilon L^{-1} (vM) \cdot\nabla_x \rho^{(2)}+O\left(\epsilon^2\right), 
     \end{equation*}
     and inserting in the macro equation \cref{macro std high stages} for $j=2$ leads to (since $a_{21}=0$) 
    \begin{equation*}
        \rho^{(2)}=\rho^n - \frac{\Tilde{a}_{21}}{a_{22}} \Delta t \left< v \otimes L^{-1}\left( (I-\Pi) v \nabla^2_x g^{(1)} \right) \right>_V  
        - a_{22} \Delta t \nabla_x \cdot \left( \left< v \otimes L^{-1} (v M) \right>_V \nabla_x  \rho^{(2)} \right) + O(\epsilon), 
     \end{equation*}
     which is not consistent with the  diffusion equation. Thus, for CK-ARS, asymptotic consistency cannot be attained if the initial data is not well-prepared. 
\end{remark}

\section{Space and velocity discretization}
\label{sec:5}
In this section, we present the spatial (for both non-staggered and staggered grids) and velocity discretization strategies that we employ in our numerical scheme.  

\subsection{Discrete velocity method}
For the velocity discretization, we will follow the discrete velocity method \cite{dom_jin}. Thus, the velocity domain is truncated as $v \in [-v_{\max},v_{\max}]$, and a uniform mesh is used $v_k = -v_{\max}+k\Delta v$, $k=1, \dots, N_v (N_v\in\mathbb{N}^\star)$ and $\Delta v=2v_{\max}/N_v$. Further,  $f(t,x,v)$ and $M(v)$ are represented as:
\begin{equation*}
    f_k(t,x):=f(t,x,v_k), \quad M_k:=M(v_k) \text{ for } k=1, \dots, N_v.  
\end{equation*}
Then, according to the definitions \eqref{def_crochet} and \eqref{mic-mac}, we have for  $j=1, \dots, N_v$ 
\begin{gather*}
    \rho (t,x) \approx  \frac{\sum_{k=0}^{N_v-1} f_k \Delta v}{\sum_{k=0}^{N_v-1} M_k \Delta v} \; \mbox{ and } \; 
    \left(\Pi f(t,x,v)\right)_j \approx   \frac{\sum_{k=0}^{N_v-1} f_k \Delta v}{\sum_{k=0}^{N_v-1} M_k \Delta v} M_j. 
\end{gather*}
For the presentation, we will skip the velocity part to focus on space discretization. 

\subsection{Space discretization using staggered grid}
\label{Subsec: Stag grid}
First, we will consider staggered grid to approximate 
$g^{(j)}$ and $\rho^{(j)}$ in space following \cite{doi:10.1137/07069479X}: the 
two meshes of the space interval $[0,1]$ are 
$x_i=i\Delta x$ and $x_{i+1/2} = (i+1/2)\Delta x$ for $i=0, \dots, N_x (N_x\in \mathbb{N}^\star)$, with $\Delta x=L/N_x$. Periodic boundary conditions will be considered in this section. 

The expressions for $g^{(j)}$ and $\rho^{(j)}$ in \eqref{Exp for g_std high}-\eqref{Exp for rho_std high} are spatially discretised by considering staggered grid: $\rho^{(j)}$ is stored at $x_i$ 
($\rho^{(j)}_i\approx \rho^{(j)}(x_i)$), and $g^{(j)}$ is stored at $x_{i+1/2}$ ($g^{(j)}_{i+1/2}(v)\approx g^{(j)}(x_{i+1/2}, v)$). The term $v\cdot\nabla_x g^{(k)}$ in \eqref{Exp for g_std high} and \eqref{Exp for rho_std high} is discretised in an upwind fashion as $v\cdot\nabla_x \approx v^+\cdot\mathbf{G}^-_{\textsf{upw}} + v^-\cdot\mathbf{G}^+_{\textsf{upw}}$ where $v^\pm=(v\pm |v|)/2$, $\mathbf{G}^\pm_{\textsf{upw}}$ denote the $N_x\times N_x$ matrices that approximate $\nabla_x$. For instance, the first order version is 
\begin{equation}
\label{S Space nabla v+ upwind}
   \mathbf{G}^-_{\textsf{upw}} = \frac{1}{\Delta x} \mathsf{circ}([-1,\underline{1}]), \;\;\;\;\;\;\;\;\;   
    \mathbf{G}^+_{\textsf{upw}} = \frac{1}{\Delta x} \mathsf{circ}([\underline{-1},1]), 
\end{equation}
where the notation $\mathsf{circ}$ is defined in \cref{Sec: App Matrix notation}. 
With these notations, we get 
\begin{equation*}
    \left(v \partial_x g^{(j)} \right)_{x_{i+1/2}} \approx v^+ \frac{g^{(j)}_{i+\frac{1}{2}}-g^{(j)}_{i-\frac{1}{2}}}{\Delta x} + v^- \frac{g^{(j)}_{i+\frac{3}{2}}-g^{(j)}_{i+\frac{1}{2}}}{\Delta x} = \left( \left( v^+\mathbf{G}^-_{\textsf{upw}} + v^-\mathbf{G}^+_{\textsf{upw}} \right) g^{(j)} \right)_{i}, 
\end{equation*}
where in the last term, the $i$ index has to be understood as the $i$-th component of the vector.  
Instead of first order upwind discretization, one can also use high order upwind discretizations so that the matrices $\mathbf{G}^\pm_{\textsf{upw}}$ will be different. Further, the term $vM\cdot \nabla_x \rho^{(k)}$ in \eqref{Exp for g_std high}-\eqref{Exp for rho_std high} and the terms of the form $\nabla_x \cdot \left< \left( \cdot \right)\right>_V$ in \eqref{Exp for rho_std high} are discretised using second order central differences as in \cite{doi:10.1137/07069479X}. 
In particular, the term $vM\cdot \nabla_x \rho^{(k)}$ is approximated by 
\begin{equation}
\label{S Space nabla central g}
    \left(\! vM  \partial_x \rho^{(k)} \!\right)_{x_{i+1/2}} \!\!\!\!\!\approx vM \frac{\rho^{(k)}_{i+1}\!-\!\rho^{(k)}_{i}}{\Delta x} \!=\! 
    \left(\! vM \mathbf{G}_{\textsf{cen}_g} \rho^{(k)}\! \right)_{i}, \;  \mathbf{G}_{\textsf{cen}_g} \!\!=\! \frac{1}{\Delta x} \mathsf{circ}([\underline{-1},1]). 
\end{equation}
Finally, the gradient terms $\nabla_x \cdot \left< \left( \cdot \right)\right>_V$ in \eqref{Exp for rho_std high} are approximated as follows 
\begin{equation}
\label{S Space nabla central rho}
    \left( \partial_x \left<  \cdot \right>_V \right)_{x_{i}}\!\!\! =\!\! \frac{ \left(\left<  \cdot \right>_V\right)_{i+1/2}-\left( \left<  \cdot \right>_V\right)_{i-1/2}}{\Delta x} \!=\!\! \left( \mathbf{G}_{\textsf{cen}_{\rho}} \left<  \cdot \right>_V \right)_{{i}}, \mathbf{G}_{\textsf{cen}_{\rho}} \!\!=\!\! \frac{1}{\Delta x} \mathsf{circ}([-1,\underline{1}]). 
\end{equation}
Again, high order centered finite differences methods can be used so that it will give different expressions 
for $\mathbf{G}_{\textsf{cen}_{\rho}}$ and $\mathbf{G}_{\textsf{cen}_{g}}$. 
Let us remark that the term $\nabla_x \cdot \nabla_x = \nabla^2_x$ in \eqref{Exp for rho_std high} is approximated  by $\mathbf{G}_{\textsf{cen}_{\rho}} \mathbf{G}_{\textsf{cen}_g}$, ie  $\mathbf{G}_{\textsf{cen}_{\rho}} \mathbf{G}_{\textsf{cen}_g}=\frac{1}{\Delta x^2} \mathsf{circ}([1,\underline{-2},1])$, which gives the standard  second order approximation of the Laplacian. \\
To ease the reading, we present the fully discrete scheme for first order ARS$(1,1,1)$ but the generalization to high order can be done using the elements of Section \ref{sec:3} \begin{eqnarray*}
g^{n+1}\!&=&\! \left( \epsilon^2 I - \Delta t L \right)^{-1}\! \left( \epsilon^2 g^n \!-\! \epsilon \Delta t \left( I-\Pi \right) \left( v^+\mathbf{G}^-_{\textsf{upw}} \!+\! v^-\mathbf{G}^+_{\textsf{upw}} \right) g^n \!-\! \epsilon \Delta t  vM  \mathbf{G}_{\textsf{cen}_g} \rho^{n+1} \right)\\
\rho^{n+1} &=& \left( I -\Delta t^2 \mathbf{G}_{{\textsf{cen}}_{\rho}} \left( \left< v \otimes \left( \epsilon^2 I - \Delta t L \right)^{-1} \left(vM\right) \right>_V \mathbf{G}_{{\textsf{cen}}_g} \right) \right)^{-1} \times \\
&&\hspace{-0.6cm} \left( \rho^n \!-\!\Delta t \mathbf{G}_{{\textsf{cen}}_{\rho}} \left< v \left( \epsilon^2 I - \Delta t L \right)^{-1} \! \left( \epsilon g^n \!-\! \Delta t \left( I-\Pi \right) \left( \left( v^+\mathbf{G}^-_{\textsf{upw}} \!+\! v^-\mathbf{G}^+_{\textsf{upw}} \right) g^n \right)  \right)\right>_V  \right).  
\end{eqnarray*}

\subsection{Space discretization using non-staggered grid}
\label{Subsec: Non-stag grid}
We also address the case of non-staggered grids 
which may be more appropriate when high dimensions are 
considered in space since only one spatial mesh is used: $x_i = i \Delta x$, for $i=0,1,..,N_x$, with $\Delta x = L/N_x$. 
Let $g^{(j)}$ and $\rho^{(j)}$ in \eqref{Exp for g_std high}-\eqref{Exp for rho_std high} $\forall j \in \{1,2,..,s \}$ be approximated in space by $g^{(j)}_i(v)\approx g^{(j)}(x_i,v)$ and $\rho^{(j)}_i\approx \rho^{(j)}(x_i)$. The term $v\cdot\nabla_x g^{(k)}$ in \eqref{Exp for g_std high}-\eqref{Exp for rho_std high} is discretised in an upwind fashion as $v\cdot\nabla_x = v^+\mathbf{G}^-_{\textsf{upw}} + v^-\mathbf{G}^+_{\textsf{upw}}$, where $v^\pm=({v\pm|v|})/2$. Here,  $\mathbf{G}^\pm_{\textsf{upw}}$ denote the matrices that represent an upwind approximation 
of $\nabla_x$. For instance, the definition \eqref{S Space nabla v+ upwind} can be used, but also its third order version 
\begin{equation}
\label{NS Space nabla v+ upwind}
   \mathbf{G}^-_{\textsf{upw}} = \frac{1}{6\Delta x} \mathsf{circ}([1,-6,\underline{3},2]), \;\;\;\;\;\;\;\;\;\;
    \mathbf{G}^+_{\textsf{upw}} = \frac{1}{6\Delta x} \mathsf{circ}([-2,\underline{-3},6,-1]),  
\end{equation}
where $\mathsf{circ}$ represents the matrix notation described in \cref{Sec: App Matrix notation} can be used. The term $vM\cdot \nabla_x \rho^{(k)}$ in \eqref{Exp for g_std high}-\eqref{Exp for rho_std high} and the terms of the form $\nabla_x \cdot \left< \left( \cdot \right)\right>_V$ in \eqref{Exp for rho_std high} are discretised in central fashion, since these terms act as source in \cref{Exp for g_std high} and diffusion in \eqref{Exp for rho_std high}. Here, $\nabla_x$ is approximated by  central differences as in \eqref{S Space nabla central rho} or \eqref{S Space nabla central g} but in the non-staggered case, the same  matrix can be used for both terms. 
As an example, the fourth order central difference produces: 
\begin{equation}
\label{NS Space nabla central}
    \mathbf{G}_{\textsf{cen}} = \frac{1}{12 \Delta x} \mathsf{circ}([1,-8,\underline{0},8,-1]). 
\end{equation}
The term $\nabla_x \cdot \nabla_x = \nabla^2_x$ in \eqref{Exp for rho_std high} is discretised as the matrices product $\mathbf{G}_{\textsf{cen}}^2=\mathbf{G}_{\textsf{cen}} \mathbf{G}_{\textsf{cen}}$. 
Like in the staggered grid case, we present the fully discrete scheme for first order ARS$(1,1,1)$ time discretization to ease the reading:
\begin{eqnarray*}
g^{n+1}\!&=&\! \left( \epsilon^2 I - \Delta t L \right)^{-1}\! \left( \epsilon^2 g^n \!-\! \epsilon \Delta t \left( I-\Pi \right) \left( v^+\mathbf{G}^-_{\textsf{upw}} + v^-\mathbf{G}^+_{\textsf{upw}} \right) g^n\right. \left. \!\!- \epsilon \Delta t  vM  \mathbf{G}_{\textsf{cen}} \rho^{n+1} \right)\\
\rho^{n+1} &=& \left( I -\Delta t^2 \mathbf{G}_{\textsf{cen}} \left( \left< v \otimes \left( \epsilon^2 I - \Delta t L \right)^{-1} \left(vM\right) \right>_V \mathbf{G}_{\textsf{cen}} \right) \right)^{-1}\times \\
&&\hspace{-0.6cm} \left( \rho^n -\Delta t \mathbf{G}_{\textsf{cen}} \left< v \left( \epsilon^2 I - \Delta t L \right)^{-1}  \left( \epsilon g^n - \Delta t \left( I-\Pi \right) \left( \left( v^+\mathbf{G}^-_{\textsf{upw}} + v^-\mathbf{G}^+_{\textsf{upw}} \right) g^n \right)  \right)\right>_V  \right) 
\end{eqnarray*}

\begin{remark}
    We know that the term $\sum_{k=1}^{j} a_{jk} \frac{\Delta t}{\epsilon} \nabla_x \cdot \left< v g^{(k)} \right>_V$ in \eqref{macro std high stages} is split into first $j-1$ and last $j$ contributions, and $g^{(j)}$ is substituted for the last $j$ contribution, as in \eqref{Exp for rho_std high}. The gradient in $\sum_{k=1}^{j-1} a_{jk} \frac{\Delta t}{\epsilon} \nabla_x \cdot \left< v g^{(k)} \right>_V$ of \eqref{Exp for rho_std high} is discretised using $\mathbf{G}_{\textsf{cen}_{\rho}}$. Further, the substitution of $g^{(j)}$ for the last $j$ hints the combination of $\nabla_x \cdot \nabla_x$ as $\nabla^2_x$ for the terms of $g^{(j)}$ involving $\nabla_x g$ and $\nabla_x \rho$. However, if we choose a spatial discretization for $\nabla^2_x$ as $\mathbf{G}_{\textsf{diff}}$, then these terms will experience $\mathbf{G}_{\textsf{cen}_{\rho}} \mathbf{G}_{\textsf{cen}_g}$ for the first $j-1$ contributions and $\mathbf{G}_{\textsf{diff}}$ for the last $j$ contribution of the $\rho^{(j)}$ update equation. This disrupts the ODE structure present in RK time discretization, and hence reduction to first order time accuracy was observed numerically. Therefore, in order to retain high order time accuracy, it is important to carry out the space discretization carefully. Hence, we do not introduce a different discretization for $\nabla^2_x$, and we retain $\mathbf{G}_{\textsf{cen}_{\rho}}  \mathbf{G}_{\textsf{cen}_g}$ even for the last $j$ contribution of $\rho^{(j)}$ equation.  
\end{remark}
\begin{remark}
\label{Space g=O(epsilon)}
The matrices introduced for spatial discretization do not change the Chapman-Enskog expansion so that the AP property is still true in the fully discrete form. Thus, we have $g^{(k)}=\epsilon L^{-1}(vM) \mathbf{G}_{\textsf{cen}_g} \rho^{(k)} + O(\epsilon^2)$ for $k \in \{1,\dots,s\}$ by using type A. For CK-ARS with well-prepared data, we have $g^{(k)}=\epsilon L^{-1}(vM) \mathbf{G}_{\textsf{cen}_g} \rho^{(k)} + O(\epsilon^2)$ for $k \in \{2,\dots,s\}$.
    Inserting this in macro equation, we get the corresponding RK scheme for the diffusion 
    \begin{gather*}
        \rho^{(j)}=\rho^n - \Delta t \sum_{k=1}^j a_{jk} \mathbf{G}_{\textsf{cen}_{\rho}} \left( \left< v \otimes L^{-1}(vM) \right>_V \mathbf{G}_{\textsf{cen}_g} \rho^{(k)} \right) + O(\epsilon).
    \end{gather*}
\end{remark}


\section{Extensions to other collision operator and inflow boundary problems}
\label{sec:6}
In this section, we show that our high order AP schemes can be extended to other problems involving advection-diffusion asymptotics and inflow boundaries. 

\subsection{Advection-diffusion asymptotics}
In this part, an advection-diffusion collision operator is considered (see 
\cite{jin_pareschi_toscani, JANG2015199}),
\begin{equation} 
\label{L: Adv-Diff} 
\mathcal{L}f := Lf +\epsilon vM \cdot A    \left< f \right>_V , \quad A\in \mathbb{R}^d, \;\; |\epsilon A|<1, 
\end{equation}
where $L$ denotes a collision satisfying the properties 
listed in Section \ref{sec:2}. A famous simple example is  $Lf=\left< f \right>_V M -f$. 

Using the notations introduced in Section \ref{sec:2}, 
we can derive the micro-macro model satisfied by $\rho=\left< f \right>_V$ and $g=f-\rho M$ by applying  $\Pi$ and $I-\Pi$ to \cref{BGK eqn} with collision $\mathcal{L}$ to get the macro and micro equations in this context 
\begin{eqnarray}
\label{macro std L}
    \partial_t \rho + \frac{1}{\epsilon} \nabla_x \cdot  \left< vg\right>_V &=& 0, \\
\label{micro std L}
    \partial_t g + \frac{1}{\epsilon} \left( I-\Pi \right) \left( v\cdot\nabla_x g \right) + \frac{1}{\epsilon}vM\cdot\nabla_x\rho &=& \frac{1}{\epsilon^2} L g + \frac{1}{\epsilon} v M \cdot A \rho.
\end{eqnarray}
A Chapman-Enskog expansion can be performed to get 
$g = \epsilon L^{-1}(vM) \cdot \nabla_x\rho - \epsilon L^{-1}(vM) \cdot A \rho + {\cal O}(\epsilon^2)$.
Inserting this in the macro equation \eqref{macro std L} enables to obtain an advection-diffusion equation in the limit $\epsilon \rightarrow 0$:
\begin{equation}
\label{L: limit eqn}
    \partial_t \rho +  \nabla_x \cdot \left( \left< v \otimes L^{-1} (vM) \right>_V \nabla_x  \rho \right)    - \nabla_x \cdot \left( \left< v \otimes L^{-1} (vM) \right>_V A  \rho \right)  = 0.
\end{equation}
The goal is to design a uniformly stable high order time integrators for \eqref{macro std L}-\eqref{micro std L} so that they degenerate into a high order time integrator 
for \eqref{L: limit eqn} as $\epsilon\to 0$. 
The extension of the schemes introduced in Section \ref{sec:3} will lead to an IMEX discretization of the asymptotic model \eqref{L: limit eqn}, where 
the advection term is treated explicitely while  
the diffusion term is implicit.

\subsubsection{High order time integrator}
In this subsection, we present the discretization of macro and micro equations \eqref{macro std L}-\eqref{micro std L}. As in Section \ref{sec:3}, in the micro equation, we treat $\frac{1}{\epsilon^2}Lg$ implicitly to ensure uniform stability and the additional term $\frac{1}{\epsilon} v M \cdot A \rho $ explicitly since it will be stabilized by the implicit treatment of the stiffest term. Regarding the macro equation and the remaining terms in micro equation, we follow the lines from previous Section \ref{sec:3}. We thus obtain the following high order IMEX RK scheme to approximate 
\eqref{macro std L}-\eqref{micro std L} 
\begin{gather}
\label{macro std L stages}
    \rho^{(j)}=\rho^n - \sum_{k=1}^j a_{jk} \frac{\Delta t}{\epsilon} \nabla_x \cdot  \left< v g^{(k)}\right>_V, \\
\label{micro std L stages}
    g^{(j)} \!\!=\! g^n \!\!-\! \!\frac{\Delta t}{\epsilon}\! \Big[\!\,\sum_{k=1}^{j-1} \Tilde{a}_{jk} \mathcal{T}\!g^{(k)} \!\!+\! \sum_{k=1}^j \!a_{jk}  vM\!\!\cdot\!\!\nabla_x\rho^{(k)}\!\! -\!\! \sum_{k=1}^j \! \frac{a_{jk}}{\epsilon} Lg^{(k)}\! -\!\! \sum_{k=1}^{j-1} \Tilde{a}_{jk}  v M\!\! \cdot \!\! A \rho^{(k)}\!\Big], 
\end{gather}
where the coefficients $a_{jk}, \tilde{a}_{jk}$ are given by the Butcher tableaux. As in Section \ref{sec:3}, some calculations are required to make the algorithm explicit. First, we have 
\begin{equation}    
\label{L g_j}
g^{(j)}\!=  \mathcal{I}^{(j)}_{\epsilon, \Delta t} \left( \epsilon^2 g^{n} \!\! -\! \epsilon \Delta t \Big[\sum_{k=1}^{j-1} \Tilde{a}_{jk} \mathcal{T}g^{(k)} \!\!+\!\sum_{k=1}^{j} a_{jk} vM \cdot \nabla_x \rho^{(k)} - \frac{1}{\epsilon} \sum_{k=1}^{j-1} a_{jk} L g^{(k)} \!-\sum_{k=1}^{j-1} \Tilde{a}_{jk} v M \cdot A \rho^{(k)}\Big] \right), 
\end{equation}
with  
$\mathcal{T}\!g^{(k)}\!\!=\!\left( I\!-\!\Pi \right) \left( v\!\cdot\!\nabla_x g^{(k)} \right)$ and $\mathcal{I}^{(j)}_{\epsilon, \Delta t}\!\!=\!\left( \epsilon^2 I \!-\! a_{jj} \Delta t L \right)^{-1}$.  
Then, $\rho^{(j)}$ is obtained by inserting $g^{(j)}$  given by \eqref{L g_j} in the macro equation \eqref{macro std L stages} to get 
\begin{eqnarray}
\label{L rho_j}
    \ \ \rho^{(j)} &=& \left( I -a_{jj}^2 \Delta t^2  \nabla_x \cdot \left(\mathcal{D}^{(j)}_{\epsilon, \Delta t} \nabla_x \right) \right)^{-1} \left( \rho^n - \sum_{k=1}^{j-1} a_{jk}\frac{\Delta t}{\epsilon} \nabla_x \cdot  \left< v g^{(k)}\right>_V \right. \\
    && \left.  - a_{jj} \Delta t \nabla_x \cdot \left< v \mathcal{I}^{(j)}_{\epsilon, \Delta t} \left( \epsilon g^n - \sum_{k=1}^{j-1} \Tilde{a}_{jk} \Delta t \mathcal{T}g^{(k)} - \sum_{k=1}^{j-1} a_{jk} \Delta t  vM\cdot\nabla_x\rho^{(k)} \right. \right. \right. \nonumber \\ 
    && \left. \left. \left. + \frac{1}{\epsilon} \sum_{k=1}^{j-1} a_{jk} \Delta t L g^{(k)} + 
    \sum_{k=1}^{j-1} \Tilde{a}_{jk} \Delta t v M \cdot A \rho^{(k)} \right) \right>_V  \right), \nonumber 
\end{eqnarray}
where $\mathcal{D}^{(j)}_{\epsilon, \Delta t} = \langle v \otimes \left( \epsilon^2 I - a_{jj} \Delta t L \right)^{-1} \left(vM\right) \rangle_V$. Thus,  
$\rho^{(j)}$ can be updated by using \eqref{L rho_j}  and $g^{(j)}$ can be found explicitly by using 
\eqref{L g_j}. 

\subsubsection{Asymptotic preserving property}
This part is dedicated to the asymptotic preserving property of the scheme \eqref{L rho_j}-\eqref{L g_j}. We first show the AP property of type A time integrator, and we later remark how this property is true for the CK-ARS time integrator with well-prepared initial data. First we have  
\begin{lemma}
\label{Lem AP g_j L}
    If $g^n=O(1)$ and $g^{(k)}=O(\epsilon), \forall k \in \{1,2, \dots,j-1 \}$, then $g^{(j)}=O(\epsilon), \forall j \in \{ 2,3,..,s\}$ for small $\epsilon$. In particular, we have $\forall j\in \{ 2,3,..,s\}$ 
    \begin{equation}
    \label{AP g_j L}
        g^{(j)}=\epsilon \sum_{k=1}^{j} \frac{a_{jk}}{a_{jj}} L^{-1} (vM) \cdot \nabla_x \rho^{(k)} - \sum_{k=1}^{j-1} \frac{a_{jk}}{a_{jj}} g^{(k)} - \epsilon \sum_{k=1}^{j-1} \frac{\Tilde{a}_{jk}}{a_{jj}} L^{-1} (vM) \cdot A \rho^{(k)} + O(\epsilon^2). 
    \end{equation}
\end{lemma}
\begin{proof}
    Plugging in \cref{L g_j} the expansion~\eqref{exp_Ij} of  $\mathcal{I}^{(j)}_{\epsilon, \Delta t}$ given by \cref{def_not_I}, 
     along with the assumptions stated in the Lemma, we obtain \eqref{AP g_j L} from which we deduce $g^{(j)}={\cal O}(\epsilon)$ 
    for all $j\in \{ 2,3,..,s\}$. 
\end{proof}
\begin{remark}
    For type A time integrator,  if $g^n={\cal O}(1)$, we have from $\eqref{L g_j}$:
        \begin{equation*}
            g^{(1)}=\epsilon \frac{a_{11}}{a_{11}} vM \cdot \nabla_x \rho^{(1)} + O(\epsilon^2)=O(\epsilon). 
        \end{equation*}
        This satisfies the induction hypothesis in Lemma \ref{Lem AP g_j L}. Further, \cref{AP g_j L} holds by omitting $\sum_{k=1}^{j-1}$ terms for $j=1$. Thus, \cref{AP g_j L} is true for $j \in \{1,2,..,s \}$. 
\end{remark}
Lemma \ref{Lem AP g_j L} enables to get an expansion of $g^{(j)}$ that can be inserted in \eqref{L rho_j} to 
identify the time discretization of the asymptotic limit. 
However, this leads to quite involved calculations 
which requires to introduce some notations.  
\begin{definition}
\label{Def Pi L}
    For $j \in \{1,2,..,s\}$ and $k_1, m \in \{1,2,..,j\}$ we define
    \begin{equation}
    \label{Def Pi eqn L}
        \mathsf{\Pi}_{j,k_1}^m = \left< v \frac{a_{jk_1}}{a_{k_1k_1}} \left( \mathcal{S}^{k_0} \mathcal{S}^{k_1} \mathcal{S}^{k_2} \dots \mathcal{S}^{k_{m-1}}  \right) \left( \mathcal{R}^{k_m} \right) \right>_V, 
    \end{equation}
    with  
    \begin{gather*}
        \mathcal{S}^{k_0}=1, \;\;\;\;\;\;  
        \mathcal{S}^{k_l}=\sum_{k_{l+1}=1}^{k_l-1} \frac{a_{k_l k_{l+1}}}{a_{k_{l+1}k_{l+1}}} \text{ for }  l \in \{1,2,..,m-1 \},\ m\geq 2, \\
        \mathcal{R}^{k_m}=\sum_{k_{m+1}=1}^{k_m} a_{k_mk_{m+1}} L^{-1} (vM) \cdot \nabla_x \rho^{(k_{m+1})} - \sum_{k_{m+1}=1}^{k_m-1} \Tilde{a}_{k_mk_{m+1}} L^{-1} (vM) \cdot A \rho^{(k_{m+1})}.
    \end{gather*} 
    As usual, we will use the convention $\sum_{j=1}^q \equiv 0$ if $q\in \mathbb{Z}\backslash\mathbb{N}$. 
\end{definition}
The term $\mathsf{\Pi}_{j,k_1}^m$ will be useful 
in the following study and deserves some remarks: the index $m$ denotes the depth of the embedded sums, $j$ corresponds to the current stage and $k_1$ corresponds to the indexing over previous stages. We continue with the following lemma which gives an induction relation on $\mathsf{\Pi}_{j,k_1}^m$. 
\begin{lemma}
\label{Lem Pi rec rel}
    For $j \geq 2$, we have 
\begin{equation*}
        \mathsf{\Pi}_{j,j}^m = \sum_{k_1=1}^{j-1} \mathsf{\Pi}_{j,k_1}^{m-1} \text{ for } m \in \{2,3,..,j\},  \mbox{ and }
\mathsf{\Pi}_{j,k_1}^j = 0 \text{ for } k_1 \in \{1,2,..,j-1 \}. 
\end{equation*}
\end{lemma}
\begin{proof}
For the first relation, considering $k_1=j$ (with $j\geq 2$) in \eqref{Def Pi eqn L} leads to 
\begin{equation*}
        \mathsf{\Pi}_{j,j}^m = 
        \left< v \left( \mathcal{S}^{k_0} \mathcal{S}^{j} \mathcal{S}^{k_2} \dots \mathcal{S}^{k_{m-1}}  \right) \left( \mathcal{R}^{k_m} \right) \right>_V, 
\end{equation*}
    since $a_{jj} \neq 0$. Further, since $\displaystyle \mathcal{S}^{k_1=j} = \sum_{k_2=1}^{j-1} \frac{a_{jk_2}}{a_{k_2k_2}}$, we get 
    \begin{equation*}
        \mathsf{\Pi}_{j,j}^m = \left< v \sum_{k_2=1}^{j-1} \frac{a_{jk_2}}{a_{k_2k_2}} \left( \mathcal{S}^{k_0} \mathcal{S}^{k_2} ...\mathcal{S}^{k_{m-1}}  \right) \left( \mathcal{R}^{k_m} \right) \right>_V  
    \end{equation*}
    By employing the change of variables as $k_{\ell}\to  k_{\ell-1}$ for $\ell \in \{2,3,..,m \}$ in the right hand side of above expression, we get
    \begin{eqnarray*}
        \mathsf{\Pi}_{j,j}^m &=& \left< v \sum_{k_{1}=1}^{j-1} \frac{a_{jk_{1}}}{a_{k_{1}k_{1}}} \left( \mathcal{S}^{k_0} \mathcal{S}^{k_{1}} \dots \mathcal{S}^{k_{m-2}}  \right) \left( \mathcal{R}^{k_{m-1}} \right) \right>_V \nonumber\\ 
        &=& \sum_{k_{1}=1}^{j-1} \left< v\frac{a_{jk_{1}}}{a_{k_{1}k_{1}}} \left( \mathcal{S}^{k_0} \mathcal{S}^{k_{1}} \dots \mathcal{S}^{k_{m-2}}  \right) \left( \mathcal{R}^{k_{m-1}} \right) \right>_V = \sum_{k_{1}=1}^{j-1} \mathsf{\Pi}_{j,k_{1}}^{m-1}, 
    \end{eqnarray*}
    which proves the first identity. 

 For the second relation, considering $m=j$ in \cref{Def Pi eqn L} leads to 
    \begin{equation*}
        \mathsf{\Pi}_{j,k_1}^j = \left< v \frac{a_{jk_1}}{a_{k_1k_1}} \left( \mathcal{S}^{k_0} \mathcal{S}^{k_1} \mathcal{S}^{k_2} ...\mathcal{S}^{k_{j-1}}  \right) \left( \mathcal{R}^{k_j} \right) \right>_V
    \end{equation*}
    We first prove the relation for $j=2$. It is clear from Definition \ref{Def Pi L} that the summation in $\mathcal{S}^{k_1}$ goes from $k_2=1$ to $k_2=k_1-1$. For $k_1=1$, the summation goes to $k_2=k_1-1=0$. Thus, since $\mathcal{S}^{k_1}$ involves $\sum_{1}^0$ for $k_1=1$, it is zero according to the convention. 
    Hence $\mathsf{\Pi}_{j,k_1}^j=0$ for $k_1=1$. \\
    We now prove the relation for $j > 2$. From Definition \ref{Def Pi L}, it can be seen that the summations in $\mathcal{S}^{k_1}$ and $\mathcal{S}^{k_2}$ go from $k_2=1$ to $k_2=k_1-1$ and $k_3=1$ to $k_3=k_2-1$ respectively. Thus, the summation in $\mathcal{S}^{k_2}$ can go to atmost $k_3=k_2-1=(k_1-1)-1=k_1-2$. Proceeding in this manner, we see that the summation in $\mathcal{S}^{k_{j-1}}$ can go to atmost $k_j=k_1-(j-1)$. \\
    For $k_1 \in \{1,2,..,j-1 \}$, $k_j=k_1-(j-1) \in \mathbb{Z}\backslash\mathbb{N}$  so that $\mathcal{S}^{k_{j-1}}=0$ 
    and hence $\mathsf{\Pi}_{j,k_1}^j=0$ for $k_1 \in \{1,2,..,j-1 \}$ which ends the proof. 
\end{proof} 
Now, we can use the previous Lemma to identify 
the asymptotic numerical scheme.  
\begin{lemma}
\label{Lem Pi AP}
    When $\epsilon \rightarrow 0$, the numerical scheme   \eqref{macro std L stages}-\eqref{micro std L stages}  degenerates into  
    \begin{equation}
    \label{rho_j Pi L}
         \rho^{(j)} = \rho^n + \Delta t \sum_{k_1=1}^j \nabla_x \cdot \left( \sum_{\ell=1}^j (-1)^\ell \mathsf{\Pi}_{j,k_1}^\ell \right) \;\;\; \text{ for } j \in \{1,2,..,s \}, 
     \end{equation}
where $\mathsf{\Pi}_{j,k_1}^\ell$ is given by \cref{Def Pi L}.
\end{lemma}
\begin{proof}
 We start with the macro equation in \cref{macro std L stages} 
    \begin{equation*}
        \rho^{(j)}=\rho^n - \sum_{k_1=1}^j a_{jk_1} \frac{\Delta t}{\epsilon} \nabla_x \cdot  \langle v g^{(k_1)}\rangle_V, 
    \end{equation*}
in which we insert $g^{(k_1)}$ given by \cref{AP g_j L} to  get 
    \begin{eqnarray*}
       \rho^{(j)}&=&\rho^n - \Delta t \sum_{k_1=1}^j  \nabla_x \cdot \left< v \frac{a_{jk_1}}{a_{k_1k_1}} \left( \sum_{k_2=1}^{k_1} a_{k_1k_2} L^{-1} (vM) \cdot \nabla_x \rho^{(k_2)} \right. \right.  \left. \left. - \sum_{k_2=1}^{k_1-1} \Tilde{a}_{k_1k_2} L^{-1} (vM) \cdot A \rho^{(k_2)} \right) \right>_V  \\ 
        &&+ \frac{\Delta t}{\epsilon} \sum_{k_1=1}^j  \nabla_x \cdot \left< v \frac{a_{jk_1}}{a_{k_1k_1}} \left( \sum_{k_2=1}^{k_1-1} a_{k_1k_2}  g^{(k_2)} \right) \right>_V + O(\epsilon) \nonumber\\
       &=& \rho^n - \Delta t \sum_{k_1=1}^j  \nabla_x \cdot  \left< v \frac{a_{jk_1}}{a_{k_1k_1}} \left( \mathcal{S}^{k_0} \mathcal{R}^{k_1} \right) \right>_V  + \frac{\Delta t}{\epsilon} \sum_{k_1=1}^j  \nabla_x \cdot \left< v \frac{a_{jk_1}}{a_{k_1k_1}} \left( \sum_{k_2=1}^{k_1-1} a_{k_1k_2}  g^{(k_2)} \right) \right>_V + O(\epsilon)  \\
        &=& \rho^n - \Delta t \sum_{k_1=1}^j  \nabla_x \cdot  \mathsf{\Pi}_{j,k_1}^1  + \frac{\Delta t}{\epsilon} \sum_{k_1=1}^j  \nabla_x \cdot \left< v \frac{a_{jk_1}}{a_{k_1k_1}} \left( \sum_{k_2=1}^{k_1-1} a_{k_1k_2}  g^{(k_2)} \right) \right>_V + O(\epsilon).
    \end{eqnarray*}
    Inserting $g^{(k_2)}$ from \cref{AP g_j L} in the above equation and simplifying as before, we get,
    \begin{equation*}
        \rho^{(j)}=\rho^n - \Delta t \sum_{k_1=1}^j  \nabla_x \cdot \left( \mathsf{\Pi}_{j,k_1}^1  - \mathsf{\Pi}_{j,k_1}^2 \right)  - \frac{\Delta t}{\epsilon} \sum_{k_1=1}^j  \nabla_x \cdot \left< v \frac{a_{jk_1}}{a_{k_1k_1}} \left( \sum_{k_2=1}^{k_1-1} \frac{a_{k_1k_2}}{a_{k_2k_2}} \sum_{k_3=1}^{k_2-1} a_{k_2k_3} g^{(k_3)} \right) \right>_V + O(\epsilon) .
    \end{equation*}
    This procedure can be continued $(j-1)$ times to finally get,
    \begin{eqnarray*}
        \rho^{(j)} &=& \rho^n + \Delta t \sum_{k_1=1}^j  \nabla_x \cdot \left( \sum_{\ell=1}^{j-1} (-1)^\ell \mathsf{\Pi}_{j,k_1}^\ell \right) \\ 
        && - (-1)^{j-1} \frac{\Delta t}{\epsilon}  \sum_{k_1=1}^j  \nabla_x \cdot \left< v \frac{a_{jk_1}}{a_{k_1k_1}} \left( \sum_{k_2=1}^{k_1-1} \frac{a_{k_1k_2}}{a_{k_2k_2}} \dots 
        \sum_{k_{j-1}=1}^{k_{j-2}-1} \frac{a_{k_{j-2}k_{j-1}}}{a_{k_{j-1}k_{j-1}}} \sum_{k_j=1}^{k_{j-1}-1} a_{k_{j-1}k_j} g^{(k_j)} \right) \right>_V + O(\epsilon) \\
        &=& \rho^n + \Delta t \sum_{k_1=1}^j  \nabla_x \cdot \left( \sum_{\ell=1}^{j-1} (-1)^\ell \mathsf{\Pi}_{j,k_1}^\ell \right) \\
        && -(-1)^{j-1} \frac{\Delta t}{\epsilon}  \sum_{k_1=1}^j  \nabla_x \cdot \left< v \frac{a_{jk_1}}{a_{k_1k_1}} \left( \mathcal{S}^{k_0} \mathcal{S}^{k_1}  \dots \mathcal{S}^{k_{j-2}} \sum_{k_j=1}^{k_{j-1}-1} a_{k_{j-1}k_j} g^{(k_j)} \right) \right>_V + O(\epsilon). 
    \end{eqnarray*}
    We know from \cref{Def Pi L} that the summations in $\mathcal{S}^{k_1}$ and $\mathcal{S}^{k_2}$ go from $k_2=1$ to $k_2=k_1-1$ and $k_3=1$ to $k_3=k_2-1$ respectively. Thus, the summation in $\mathcal{S}^{k_2}$ can go to atmost $k_3=k_2-1=(k_1-1)-1=k_1-2$. Proceeding in this manner, we see that the summations in $\mathcal{S}^{k_{j-2}}$ and $\sum_{k_j=1}^{k_{j-1}-1} a_{k_{j-1}k_j} g^{(k_j)}$ go to atmost $k_{j-1}=k_1-(j-2)$ and $k_j=k_1-(j-1)$ respectively. \\
    Since the summation in $k_1$ goes to atmost $j$ in the above equation, $k_j$ in the term $\sum_{k_j=1}^{k_{j-1}-1} a_{k_{j-1}k_j} g^{(k_j)}$ goes to atmost $k_j=k_1-(j-1)=j-(j-1)=1$, and $k_{j-1}$ in $\mathcal{S}^{k_{j-2}}$ goes to atmost $k_{j-1}=k_1-(j-2)=j-(j-2)=2$ and so on. Thus, only $k_j=1$ remains in the last summation so that $\sum_{k_j=1}^{k_{j-1}-1} a_{k_{j-1}k_j} g^{(k_j)} =a_{21} g^{(1)} =\epsilon a_{21} L^{-1}(vM) \cdot \nabla_x \rho^{(1)} + {\cal O}(\epsilon^2) = \frac{a_{21}}{a_{11}} \epsilon a_{11} L^{-1}(vM) \cdot \nabla_x \rho^{(1)} + {\cal O}(\epsilon^2) = \epsilon \mathcal{S}^{k_{j-1}} \mathcal{R}^{k_{j}} + {\cal O}(\epsilon^2) $. 
    Thus, we have 
    \begin{eqnarray*}
        \rho^{(j)} &=& \rho^n + \Delta t \sum_{k_1=1}^j  \nabla_x \cdot \left( \sum_{\ell=1}^{j-1} (-1)^\ell \mathsf{\Pi}_{j,k_1}^\ell \right) \nonumber\\
        &&- (-1)^{j-1} \Delta t \sum_{k_1=1}^j \nabla_x \cdot \left< v \frac{a_{jk_1}}{a_{k_1k_1}} \left( \mathcal{S}^{k_0} \mathcal{S}^{k_1}  \dots \mathcal{S}^{k_{j-1}} \mathcal{R}^{k_j} \right) \right>_V + O(\epsilon)\\
        &=& \rho^n + \Delta t \sum_{k_1=1}^j  \Bigg[ \nabla_x \cdot \left( \sum_{\ell=1}^{j-1} (-1)^\ell \mathsf{\Pi}_{j,k_1}^\ell \right) + \nabla_x \cdot \left( (-1)^j \mathsf{\Pi}_{j,k_1}^j \right) \Bigg] + O(\epsilon).
    \end{eqnarray*} 
\end{proof}
We can now prove the asymptotic property of the scheme \eqref{macro std L stages}-\eqref{micro std L stages}. 
\begin{theorem}
\label{thm AP}
   When $\epsilon \rightarrow 0$, the scheme   \eqref{macro std L stages}-\eqref{micro std L stages}  degenerates into
    \begin{multline}
    \label{rho_j AP L}
         \rho^{(j)} = \rho^n - \Delta t \sum_{k=1}^j a_{jk} \nabla_x \cdot \left( \left<v \otimes L^{-1} (vM) \right>_V \nabla_x \rho^{(k)} \right)\\ 
         + \Delta t\sum_{k=1}^{j-1} \Tilde{a}_{jk} \nabla_x \cdot \left( \left<v \otimes L^{-1} (vM) \right>_V A \rho^{(k)} \right),  \mbox{ for } j \in \{1,2,\dots,s \}.  
     \end{multline}
\end{theorem}
\begin{proof}
    From Lemma \ref{Lem Pi AP}, the asymptotic limit $\epsilon \rightarrow 0$ of the macro equation in \cref{macro std L stages} is (for $j \in \{1,2,..,s \}$)
    \begin{align*}
        \rho^{(j)} &= \rho^n + \Delta t \sum_{k_1=1}^j \nabla_x \!\cdot\! \left( \sum_{\ell=1}^j (-1)^\ell \mathsf{\Pi}_{j,k_1}^\ell \right)= \rho^n + \Delta t  \nabla_x \!\cdot\! \left( \sum_{\ell=1}^j (-1)^\ell \left( \mathsf{\Pi}_{j,j}^\ell + \sum_{k_1=1}^{j-1}    \mathsf{\Pi}_{j,k_1}^\ell \right) \right) \nonumber\\
        &= \rho^n + \Delta t  \nabla_x \cdot \left( -  \mathsf{\Pi}_{j,j}^1 + \sum_{\ell=2}^j (-1)^\ell \mathsf{\Pi}_{j,j}^\ell + \sum_{\ell=1}^j (-1)^\ell \sum_{k_1=1}^{j-1}    \mathsf{\Pi}_{j,k_1}^\ell  \right).
    \end{align*}
Using the recurrence relation given by Lemma \ref{Lem Pi rec rel} and a change of indices lead to 
    \begin{align*}
        \rho^{(j)} &= \rho^n + \Delta t  \nabla_x \cdot \left( -   \mathsf{\Pi}_{j,j}^1 + \sum_{\ell=2}^j (-1)^\ell \sum_{k_1=1}^{j-1}    \mathsf{\Pi}_{j,k_1}^{\ell-1} + \sum_{\ell=1}^j (-1)^\ell \sum_{k_1=1}^{j-1}    \mathsf{\Pi}_{j,k_1}^\ell  \right)\nonumber\\
        &= \rho^n + \Delta t  \nabla_x \cdot \left( -   \mathsf{\Pi}_{j,j}^1 - \sum_{\ell=1}^{j-1} (-1)^{\ell} \sum_{k_1=1}^{j-1}    \mathsf{\Pi}_{j,k_1}^{\ell} + \sum_{\ell=1}^j (-1)^\ell \sum_{k_1=1}^{j-1}    \mathsf{\Pi}_{j,k_1}^\ell  \right)\nonumber\\
         &= \rho^n + \Delta t  \nabla_x \cdot \left( -  \mathsf{\Pi}_{j,j}^1 +  (-1)^j \sum_{k_1=1}^{j-1}    \mathsf{\Pi}_{j,k_1}^j  \right).
    \end{align*}
    From Lemma \ref{Lem Pi rec rel}, we have $\sum_{k_1=1}^{j-1} \mathsf{\Pi}_{j,k_1}^j= 0$, so that from Definition \ref{Def Pi L} we get  
    \begin{align*}
        \rho^{(j)} &= \rho^n + \Delta t  \nabla_x \cdot \left( -\mathsf{\Pi}_{j,j}^1 \right)  = \rho^n - \Delta t \nabla_x \cdot \left( \left< v \frac{a_{jj}}{a_{jj}} \mathcal{S}^{k_0}   \mathcal{R}^{k_1=j}  \right>_V \right) \nonumber\\
         &= \rho^n - \Delta t \nabla_x \cdot \left( \left< v \left( \sum_{k_{2}=1}^{k_1} a_{k_1k_{2}} L^{-1} (vM) \cdot \nabla_x \rho^{(k_{2})} \right. \right.  \left. \left. - \sum_{k_{2}=1}^{k_1-1} \Tilde{a}_{k_1k_{2}} L^{-1} (vM) \cdot A \rho^{(k_{2})} \right) \right>_V \right)_{k_1=j}\nonumber\\
       &= \rho^n - \Delta t \sum_{k_{2}=1}^{j}  a_{jk_{2}} \nabla_x \cdot \left( \left< v \otimes L^{-1} (vM) \right>_V \nabla_x \rho^{(k_{2})} \right)  +  \Delta t \sum_{k_{2}=1}^{j-1} \Tilde{a}_{jk_{2}}  \nabla_x \cdot \left( \left< v \otimes L^{-1} (vM) \right>_V A \rho^{(k_{2})} \right), 
   \end{align*}
   which ends the proof. 
\end{proof}
\begin{remark}
     For CK-ARS schemes with well-prepared initial data, we obtain $g^{(1)}=g^n=O(\epsilon)$ and $\rho^{(1)}=\rho^n$. The presentation in this section will apply for CK-ARS from the second RK stage onwards. For instance, \cref{Def Pi L} applies for CK-ARS with the following change in indexes: $j \in \{2,3,..,s \}$, $k_1,m \in \{2,3,..,j \}$ and all the summations involved start from 2 instead of 1 since $a_{11}=0$. The lemmas and theorems that follow also undergo the corresponding change in indexes, and the AP property for CK-ARS can be observed for $j \in \{2,3,..,s \}$. 
\end{remark}
\begin{remark}
    Upon incorporating the spatial matrices corresponding to staggered grid in place of the continuous gradient operator, we obtain in the limit $\epsilon \rightarrow 0$,
\begin{multline}
        \rho^{(j)} = \left( I + a_{jj}\Delta t \mathbf{G}_{\textsf{cen}_{\rho}} \left( \left< v \otimes L^{-1}(vM) \right>_V  \mathbf{G}_{\textsf{cen}_g} \right) \right)^{-1} \times\\ 
        \left( \rho^n - \sum_{k=1}^{j-1} a_{jk}\Delta t \mathbf{G}_{\textsf{cen}_{\rho}}\left( \left<v \otimes L^{-1}(vM)\right>_V \mathbf{G}_{\textsf{cen}_g} \rho^{(k)} \right) \right.  \left. + \sum_{k=1}^{j-1} \Tilde{a}_{jk} \Delta t \mathbf{G}_{\textsf{cen}_{\rho}} \left( \left< v \otimes L^{-1} (vM) \right>_V \mathbf{G}_{avg_{g}} A \rho^{(k)} \right) \right).
\end{multline}
The matrices $\mathbf{G}_{\textsf{cen}_{\rho}}, \mathbf{G}_{\textsf{cen}_g}$ are given in \cref{Subsec: Stag grid} and $\mathbf{G}_{\textsf{avg}_g}=\frac{1}{2} \mathsf{circ}([\underline{1},1])$. Thus, $A \left( \rho^{(k)} \right)_{x_{i+1/2}} = \frac{1}{2} A (\rho^{(k)}_{i+1}+\rho^{(k)}_{i}) = ( \mathbf{G}_{\textsf{avg}_g} A \rho^{(k)} )_{i}$. This results in a central discretization of the advection term in the macro equation. Thus, we obtain a consistent internal RK stage approximation of the advection-diffusion equation in the limit $\epsilon \rightarrow 0$. 
\end{remark}

\subsection{Inflow Boundaries}
So far, periodic boundary conditions were considered. In this part, we consider inflow boundary conditions for $f$ solution to \eqref{BGK eqn}
\begin{equation}
\label{inflow bc}
    f(t,x,v)=f_b(t,x,v), \quad (x,v) \in \partial\Omega \times V \text{ such that } v \cdot n(x) < 0, \quad \forall t, 
\end{equation}
where $f_b$ is a given function and $n(x)$ denotes the unitary outgoing normal vector to $\partial \Omega$.  
As mentioned in \cite{doi:10.1137/07069479X, doi:10.1137/120865513}, such boundary conditions 
cannot be adapted naturally to the standard micro-macro unknown  $\rho(t,x)$ and $g(t,x,v)$ solution to  \eqref{mic-mac}. To overcome this drawback, another  micro-macro decomposition is introduced in \cite{doi:10.1137/120865513} 
\begin{equation}
\label{mic-mac boundary}
    f=\overline{\rho}M+\overline{g}, \; 
   \overline{\rho}(t,x)=\left< f(t, x, \cdot) \right>_{V_-}\!\!, \; \left< \overline{g}(t,x,\cdot) \right>_{V_-}\!\! = 0, \;  \langle f \rangle_{V_-} \!\!= \dfrac{\int_{V_-} f d\mu}{ \int_{V_-} M d\mu}, 
\end{equation}
where the velocity domain $V_-$ is defined by 
\begin{equation}
    V_-(x)= \{ v \in V, \omega (x,v) < 0 \}, \quad V_+(x)=V \backslash V_-(x). 
\end{equation}
The function $\omega(x,v)$ extends $v \cdot n(x)$ in the interior of domain. Some examples of $\omega(x,v)$ for different geometries are provided in \cite{doi:10.1137/120865513}. It can be seen that the boundary conditions for $\overline{\rho}(t,x)$ and $\overline{g}(t,x,v)$ can be evaluated from the inflow boundary condition in \cref{inflow bc}. Indeed, for $(x,v) \in \partial \Omega \times V$ such that $v \cdot n(x) < 0$, $\forall t$, we define 
\begin{equation}
\label{mic-mac bdy rho_b g_b}
\overline{\rho}_b(t,x)= \left< f_b(t, x,\cdot) \right>_{V_-}, \;\; \overline{g}_b(t,x,v) = f_b(t,x,v) -\overline{\rho}_b(t,x) M(v).   
\end{equation}

The derivation of the micro-macro model needs to be adapted to this decomposition. The projector $\Pi^-$ is defined as $\Pi^-h=\left< h \right>_{V_-} M$. Then, substituting \cref{mic-mac boundary} into \cref{BGK eqn} and applying $\Pi^-$ and $I-\Pi^-$ enable to get the macro and micro equations:
\begin{gather}
\label{macro bdy1}
    \partial_t \overline{\rho} + \frac{1}{\epsilon} \left< vM \right>_{V_-} \cdot \nabla_x \overline{\rho} + \frac{1}{\epsilon} \nabla_x \cdot \left< v \overline{g}\right>_{V_-} = \frac{1}{\epsilon^2} \left<L\overline{g}\right>_{V_-},  \\
    \label{micro bdy}
    \partial_t \overline{g} + \frac{1}{\epsilon} \left( I-\Pi^- \right) \left( v\cdot\nabla_x \overline{g} \right) + \frac{1}{\epsilon} \left( I-\Pi^- \right) vM \cdot \nabla_x \overline{\rho} = \frac{1}{\epsilon^2} \Tilde{L}\overline{g}, 
\end{gather}
where $\Tilde{L}=\left( I-\Pi^- \right) L$. Moreover, it can be seen that $ \Tilde{L} = \left( I-\Pi^- \right) L \left( I-\Pi^- \right)=\left( I-\Pi^- \right) L \left( I-\Pi \right)$ since $\Pi^-h,\Pi h \in \mathcal{N}(L), \forall h$. \\ The macro equation \eqref{macro bdy1} turns out to be more complicated than the one obtained for standard micro-macro decomposition. It can be made simpler by using $\rho=\overline{\rho}+\left<\overline{g}\right>_V, f=\rho M - \left<\overline{g}\right>_V M + \overline{g}$, obtained from the  decompositions  \eqref{mic-mac} and \eqref{mic-mac boundary}. Applying $\Pi$ to \cref{BGK eqn} instead of $\Pi^-$, we obtain the simpler macro equation,
\begin{equation}
\label{macro bdy}
    \partial_t \rho + \frac{1}{\epsilon} \nabla_x \cdot \left< v \overline{g}\right>_{V} = 0,  
\end{equation}
and the micro-macro system that we will consider in the sequel is \eqref{micro bdy}-\eqref{macro bdy}. 

\subsubsection{Numerical scheme}
\label{subsec:boundary}
In this part, we present the fully discretized 
scheme to approximate \eqref{micro bdy}-\eqref{macro bdy}. The boundary conditions on $\overline{\rho}_b$ and $\overline{g}_b$ in \cref{mic-mac bdy rho_b g_b} will be utilised along with the relation $\rho=\overline{\rho}+\left<\overline{g}\right>_V$ that allows to link $\rho$ and $\overline{\rho}$ in the interior of the domain. We will use a staggered grid in space following \cite{doi:10.1137/120865513} and a high order scheme in time, following the strategy developed previously. To ease the reading, only the first order version will be presented. \\
First, we present the space approximation based on a staggered grid. Let us consider the space interval $[0, L]$  with two grids: $x_i\!=\!i\Delta x$ and $x_{i+1/2}=(i+1/2)\Delta x$, $\Delta x=L/(N_x-1)$. The 'interior' variables  such as $\rho, \overline{\rho}$ are stored at grid points $x_i$ with $i=1, \dots, N_x-2)$ and $\overline{g}$ is 
stored at $i+1/2=1/2, \cdots, N_x-3/2$. We also use the variable $\overline{g}_{cl}= \bar{g}\cup \bar{g}_b \in\mathbb{R}^{N_x+1}$. 
The whole domain including boundary will be considered for the micro unknown $\bar{g}$ so that the components of $\overline{g}_{cl}$  correspond to the grid indices  $i+1/2=-1/2,\cdots,N_x-1/2$. The matrices corresponding to spatial operators are given by 
\begin{eqnarray}
\label{Bupw}
\;\;\;\;\;\;\;\;  \mathbf{B}^-_{\textsf{upw}}&=& \frac{1}{\Delta x} \mathsf{circ}([\underline{-1},1])_{(N_x-1)\!\times\!(N_x+1)},  \mathbf{B}^+_{\textsf{upw}}= \frac{1}{\Delta x} \mathsf{circ}([\underline{0},-1,1])_{(N_x-1)\!\times\!(N_x+1)}, \\
\label{Bcen}
\;\;\;\;\;\;\;\;\;\;\;\;  \mathbf{B}_{\textsf{cen}_{\rho}}&=& \frac{1}{\Delta x} \mathsf{circ}([\underline{-1},1])_{(N_x-2)\times(N_x-1)}, \ 
\mathbf{B}_{\textsf{avg}}=\frac{1}{2} \mathsf{circ}([\underline{1},1])_{(N_x-2)\times(N_x-1)},\\
\label{Bceng}
\mathbf{B}_{\textsf{cen}_g}&=& \frac{1}{\Delta x} \mathsf{circ_b}([-1,\underline{1}])_{(N_x-1)\times(N_x-2)}. 
\end{eqnarray}
The $\mathsf{circ_b}$ definition is presented in \cref{Sec: App Matrix notation}. Further, we also introduce a vector containing the boundary values of $\overline{\rho}$ as $\overline{\rho}_{bd}= \frac{1}{\Delta x} \begin{bmatrix}
        -\overline{\rho}_{b_{i=0}},0,0,...,0,\overline{\rho}_{b_{i=N_x-1}}
    \end{bmatrix}^T_{(N_x-1)\times 1}$. 
We now present our scheme by using this matrix notation. For simplicity, we assume that $\overline{\rho}_{bd}$ is time invariant. 
We also use the following notations:
\begin{gather*}
    \overline{\mathcal{T}}h=\left( I-\Pi^- \right) \left( v^+ \mathbf{B}^-_{\textsf{upw}}  +   v^-\mathbf{B}^+_{\textsf{upw}} \right) h,     \overline{\mathcal{D}}_{\epsilon, \Delta t} = \bigl\langle v \bigl( \epsilon^2 I - \Delta t \Tilde{L} \bigr)^{-1} \Delta t \left( I-\Pi^- \right) (vM) \rangle_V, \\
    \overline{\mathcal{E}}_{\epsilon, \Delta t} = \bigl\langle \bigl( \epsilon^2 I - \Delta t \Tilde{L} \bigr)^{-1} \Delta t \left( I-\Pi^- \right) (vM) \bigr\rangle_V, \;\;\;
    \overline{\mathcal{I}}_{\epsilon, \Delta t} = \bigl( \epsilon^2 I - \Delta t \Tilde{L} \bigr)^{-1}, \;\;\; 
    \overline{\mathcal{J}} = \left( I-\Pi^- \right) (vM).  
\end{gather*}
The micro equation \eqref{micro bdy} is discretised in time as in the previous (periodic) case  
\begin{equation}
\label{g_boundary}
    \overline{g}^{n+1}= \overline{\mathcal{I}}_{\epsilon, \Delta t} 
    \left( \epsilon^2 \overline{g}^n 
     - \epsilon \Delta t \overline{\mathcal{T}} \overline{g}_{cl}^n -\epsilon \Delta t \overline{\mathcal{J}} \mathbf{B}_{\textsf{cen}_g} \overline{\rho}^{n+1} -\epsilon \Delta t \overline{\mathcal{J}}  \overline{\rho}_{bd} \right), 
\end{equation}
and for the macro equation \eqref{macro bdy}, we obtain 
\begin{equation*}
    \frac{\rho^{n+1}-\rho^n}{\Delta t} + \frac{1}{\epsilon} \left< v \mathbf{B}_{\textsf{cen}_{\rho}} \overline{g}^{n+1} \right>_{V} = 0
\end{equation*}
Substituting $\overline{g}^{n+1}$ in the above equation, we get
\begin{equation}
    \rho^{n+1} = \rho^n - \Delta t \mathbf{B}_{\textsf{cen}_{\rho}} \left< v  \overline{\mathcal{I}}_{\epsilon, \Delta t} \left( \epsilon \overline{g}^n  -  \Delta t \overline{\mathcal{T}} \overline{g}_{cl}^n  - \Delta t \overline{\mathcal{J}} \mathbf{B}_{\textsf{cen}_g} \overline{\rho}^{n+1} - \Delta t \overline{\mathcal{J}}  \overline{\rho}_{bd}  \right)  \right>_{V}.
\end{equation}
In index notation, we use  $\rho_i^{n+1}=\overline{\rho}_i^{n+1}+\frac{1}{2}\langle \overline{g}_{i-1/2}^{n+1} + \overline{g}_{i+1/2}^{n+1}  \rangle_V$ (since $\rho=\overline{\rho}+\langle\overline{g}\rangle_V$) to match the two grids. 
In matrix notation, this becomes $\rho^{n+1}=\overline{\rho}^{n+1}+\mathbf{B}_{\textsf{avg}} \langle \overline{g}^{n+1}\rangle_V$ with $\mathbf{B}_{\textsf{avg}}$ given by \eqref{Bcen}. Substituting this into the above equation and inserting the expression for $\overline{g}^{n+1}$ into $\mathbf{B}_{\textsf{avg}} \left< \overline{g}^{n+1}\right>_V$ enable to update the interior macro unknown 
\begin{multline}
\label{rho_boundary}
\overline{\rho}^{n+1} = \left( I - \epsilon \mathbf{B}_{\textsf{avg}} \left( \overline{\mathcal{E}}_{\epsilon, \Delta t} \mathbf{B}_{\textsf{cen}_g} \right) - \Delta t \mathbf{B}_{\textsf{cen}_{\rho}} \left( \overline{\mathcal{D}}_{\epsilon, \Delta t} \mathbf{B}_{\textsf{cen}_g} \right) \right)^{-1}\times  \\ 
\left( \rho^n - \mathbf{B}_{\textsf{avg}} \left<\overline{\mathcal{I}}_{\epsilon, \Delta t} \left( \epsilon^2 \overline{g}^n -  \epsilon \Delta t \overline{\mathcal{T}} \overline{g}_{cl}^n -\epsilon \Delta t \overline{\mathcal{J}}  \overline{\rho}_{bd}   \right)  \right>_{V}  \right. \\
\left. -\Delta t \mathbf{B}_{\textsf{cen}_{\rho}} \left< v  \overline{\mathcal{I}}_{\epsilon, \Delta t} \left( \epsilon \overline{g}^n  -  \Delta t \overline{\mathcal{T}} \overline{g}_{cl}^n  - \Delta t \overline{\mathcal{J}}  \overline{\rho}_{bd}  \right) \right>_{V} \right).  
\end{multline}
The right hand side of above expression involves only known quantities so that $\overline{\rho}^{n+1}$ can be updated from \eqref{rho_boundary} which can then be used to update  $\overline{g}^{n+1}$ in \eqref{g_boundary}. 
Then, we update $\overline{g}_{cl}^{n+1}$ thanks to the boundary conditions \eqref{mic-mac bdy rho_b g_b}, and finally $\rho^{n+1}$ can be computed from $\rho^{n+1}=\overline{\rho}^{n+1}+\mathbf{B}_{\textsf{avg}} \left< \overline{g}^{n+1}\right>_V$.  In the limit $\epsilon \rightarrow 0$, the above equation becomes, 
\begin{equation*}
    \overline{\rho}^{n+1}\!=\!\left( I + \Delta t \mathbf{B}_{\textsf{cen}_{\rho}} \left( \left< v \otimes  \Tilde{L}^{-1} \!\! \overline{\mathcal{J}} \right>_V \mathbf{B}_{\textsf{cen}_g} \right) \right)^{-1} 
    \left( \rho^n \!-\!\Delta t \mathbf{B}_{\textsf{cen}_{\rho}} \left( \left< v \otimes  \Tilde{L}^{-1} \overline{\mathcal{J}}   \right>_{V} \overline{\rho}_{bd} \right) \right) 
\end{equation*}
This is a consistent discretization of the diffusion equation in \cref{Diff eqn} since $\langle v \otimes  \Tilde{L}^{-1} \overline{\mathcal{J}}\rangle_{V} = \langle v \otimes L^{-1} (vM) \rangle_{V}=-\kappa$. Further, the high order scheme in time can be constructed in a similar manner as before. 

\section{Numerical results}
\label{sec:7}
In this section, we present the numerical validation of our high order asymptotic preserving schemes in different configurations.  
\subsection{Diffusion asymptotics}
First, we check time and space accuracy for the micro-macro scheme in the diffusion limit. 
\subsubsection{Time order of accuracy}
The spatial domain $L=[0,2\pi]$ of the problem is discretized using $N_x=50$ grid points. The velocity domain is truncated to $[-v_{\max}, v_{\max}]$ with $v_{\max}=5$ 
and we take $\Delta v=1$. 
The initial condition is: 
\begin{gather*}
    \rho(0,x)=1+\cos(x) \\
    \text{Well-prepared data (WP): } g(0,x,v) = \epsilon^2 (I-\Pi) \left( v^2 M \right) \rho (0,x) \\
    \text{Non-well prepared data (N-WP): } g(0,x,v) = (I-\Pi) \left( v^2M \right) \rho (0,x), 
\end{gather*}
with $M(v)=\frac{1}{\sqrt{2\pi}} e^{-v^2/2}$. 
Periodic boundary conditions are used on both $\rho$ and $g$. The spatial terms are discretised by using the atmost-third order accurate matrices on non-staggered grid presented in \cref{Subsec: Non-stag grid}. The final time is $T=0.5$, and the following $\Delta t$ are considered 
to validate the different high order time integrators: $\Delta t=0.5, 0.1, 0.05, 0.01, 0.005, 0.001$. The type A micro-macro schemes constructed using the Butcher tableau corresponding to DP-A$(1,2,1)$, DP2-A$(2,4,2)$ and DP1-A$(2,4,2)$ are considered. Although DP1-A$(2,4,2)$ is second order accurate, the implicit part of it when used separately is third order accurate. Further, we also consider the type CK-ARS micro-macro schemes constructed using Butcher tableau corresponding to ARS$(1,1,1)$, ARS$(2,2,2)$ and ARS$(4,4,3)$. The Butcher tableau of different time integrators utilised are presented in \cref{App: Butcher tableau}.\\ 

\begin{figure}
    \centering
    \subfloat[A N-WP, $\epsilon=1$]{\label{Anwpem0time}\includegraphics[width=0.25\textwidth]{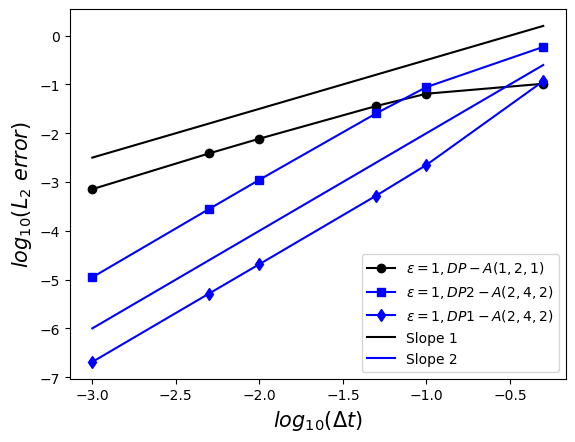}}
    \subfloat[A WP, $\epsilon=1$]{\label{Awpem0time}\includegraphics[width=0.25\textwidth]{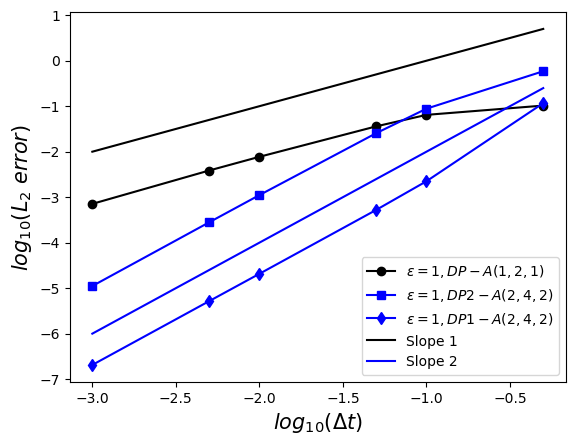}}
    \subfloat[A N-WP, $\epsilon=10^{-4}$]{\label{Anwpem4time}\includegraphics[width=0.25\textwidth]{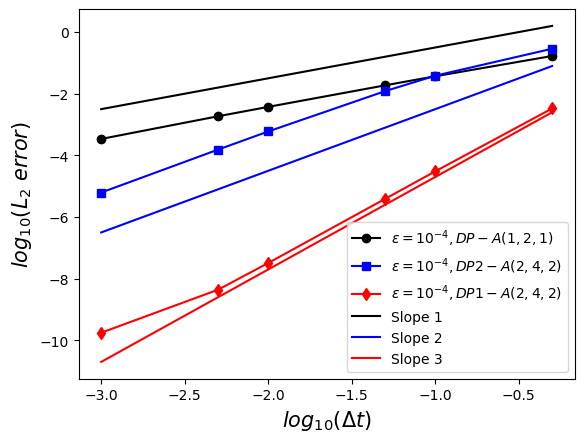}}
    \subfloat[A WP, $\epsilon=10^{-4}$]{\label{Awpem4time}\includegraphics[width=0.25\textwidth]{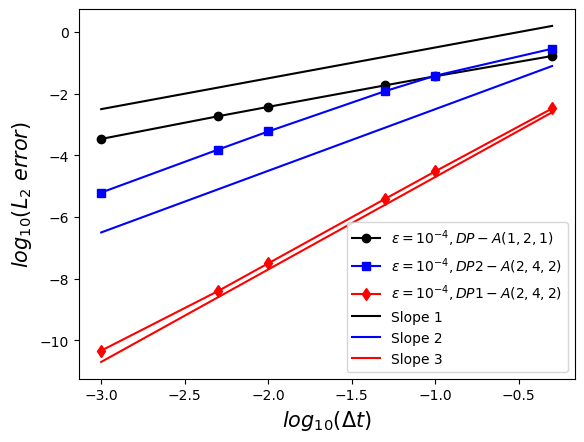}}

    \subfloat[CK N-WP, $\epsilon=1$]{\label{CKARSnwpem0time}\includegraphics[width=0.25\textwidth]{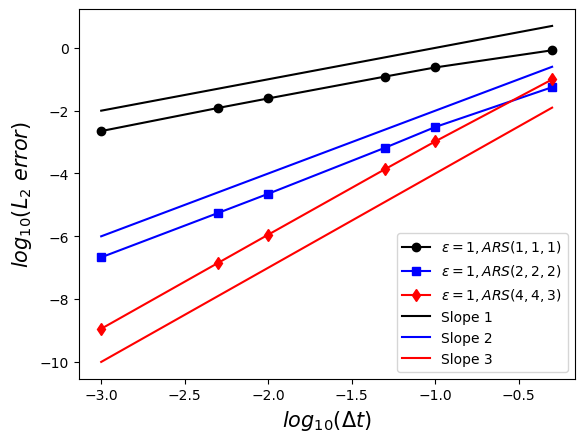}}
    \subfloat[CK WP, $\epsilon=1$]{\label{CKARSwpem0time}\includegraphics[width=0.25\textwidth]{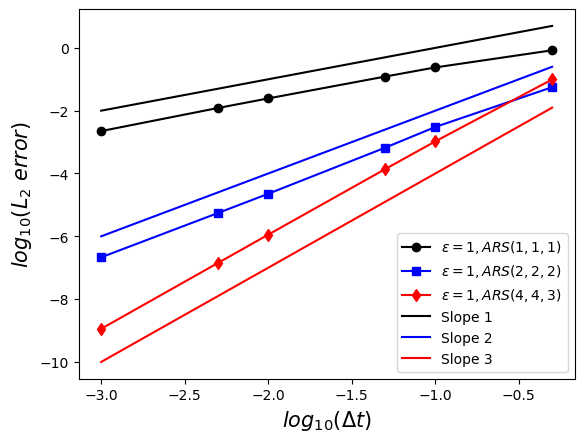}}
    \subfloat[CK N-WP, $\epsilon\!=\!10^{-4}$]{\label{CKARSnwpem4time}\includegraphics[width=0.25\textwidth]{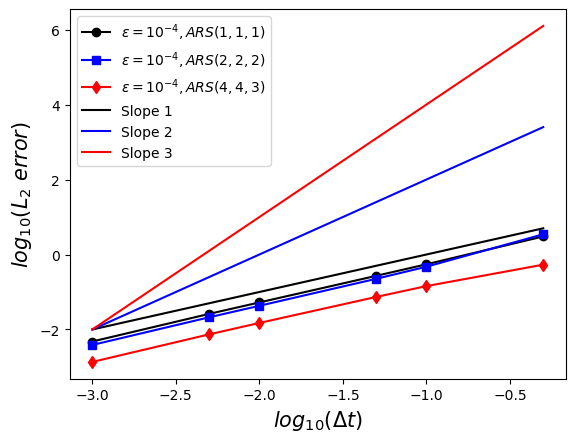}}
    \subfloat[CK WP, $\epsilon=10^{-4}$]{\label{CKARSwpem4time}\includegraphics[width=0.25\textwidth]{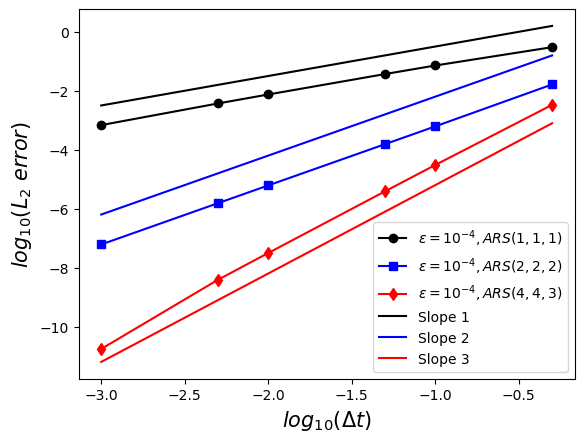}}
    \caption{Accuracy in time for different type A and CK-ARS time integrators (both WP and N-WP initial data). The reference solution is obtained from the micro-macro with $\Delta t =10^{-4}$. }
    \label{time order}
\end{figure}

\begin{figure}
    \centering
    \subfloat[A N-WP]{\label{Anwpem4timediffref}\includegraphics[width=0.33\textwidth]{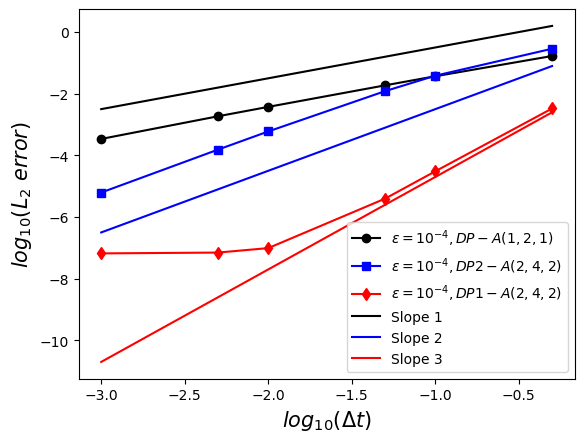}}
    \subfloat[A WP]{\label{Awpem4timediffref}\includegraphics[width=0.33\textwidth]{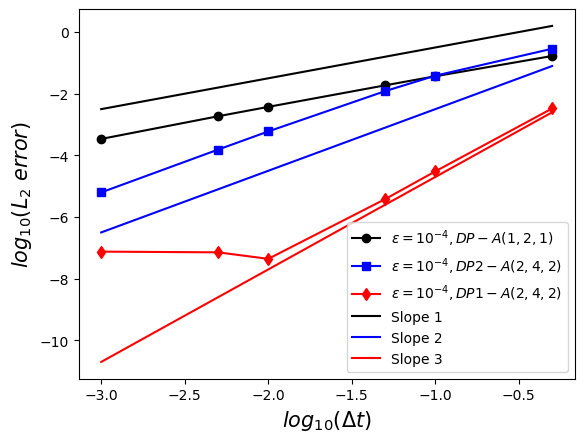}}

    \subfloat[CK N-WP]{\label{CKARSnwpem4timediffref}\includegraphics[width=0.33\textwidth]{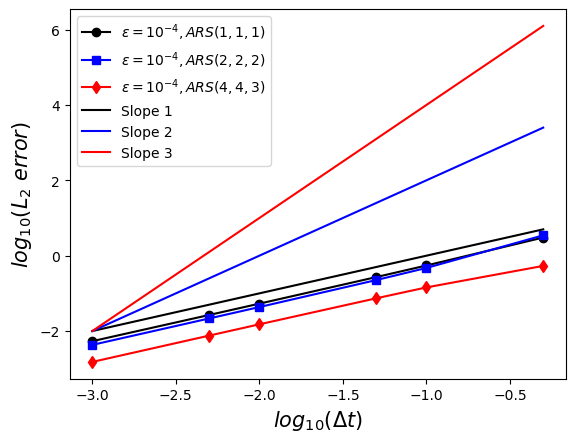}}
    \subfloat[CK WP]{\label{CKARSwpem4timediffref}\includegraphics[width=0.33\textwidth]{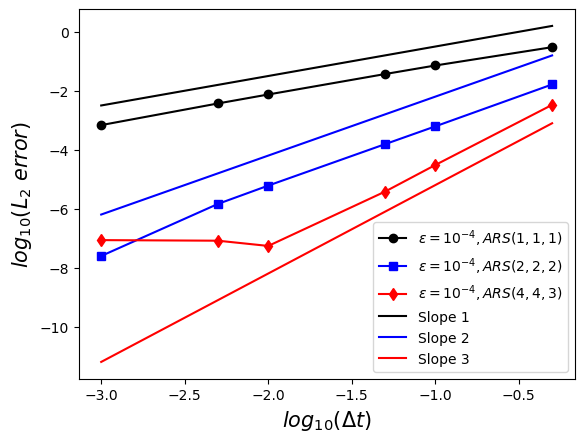}}

    \caption{Accuracy in time for different type A and CK-ARS time integrators (both WP and N-WP initial data). The reference solution is obtained from the diffusion equation with $\Delta t =10^{-4}$. }
    \label{time order diffref}
\end{figure}

In \cref{time order}, we plot the time error for the different time integrators in both WP and N-WP cases and 
for different values of $\epsilon$. Note that 
the reference solution for each curve is obtained by using the same micro-macro scheme corresponding to that curve with  $\Delta t = 10^{-4}$. 
For $\epsilon=1$, the required orders of accuracy are recovered for type A schemes with both N-WP and WP initial data, as observed in \cref{Anwpem0time,Awpem0time}. 
For $\epsilon=10^{-4}$, due to the asymptotic degeneracy of our scheme into a fully-implicit scheme for diffusion equation, only the implicit part of the Butcher tableau plays a role. Hence DP1-A$(2,4,2)$ becomes third order accurate in time, while DP-A$(1,2,1)$ and DP2-A$(2,4,2)$ are first and second order accurate respectively. This is shown in \cref{Anwpem4time,Awpem4time}. On the other hand, CK-ARS schemes with both N-WP and WP initial data for $\epsilon=1$ recover the required orders of accuracy as shown in \cref{CKARSnwpem0time,CKARSwpem0time}. However, for $\epsilon=10^{-4}$,   orders of accuracy are observed only when WP initial data are used (\cref{CKARSwpem4time}). As shown in the analyses presented in previous sections, usage of N-WP initial data for CK-ARS time integrators does not allow the asymptotic accuracy (\cref{CKARSnwpem4time}), as discussed  
in \cite{dimarco_imexrk}. \\
Since we proved the asymptotic preserving property, the diffusion solution is used as reference solution in the asymptotic regime $(\epsilon=10^{-4})$ with $\Delta t =10^{-4}$ (in \cref{time order diffref}) 
to check the orders of accuracy of high order integrators. 
The results are similar to the ones obtained for $\epsilon=10^{-4}$ in \cref{time order}, except that here we observe a plateau for third order scheme and small $\Delta t$. This is due to the ${\cal O}(\epsilon^2)$ difference between the schemes based on micro-macro and diffusion models. This error dominates ${\cal O}(\Delta t^3)$ error, and hence it is observed.    

\subsubsection{Space order of accuracy}
\begin{figure}
    \centering
    \subfloat[A N-WP]{\label{Anwpspace}\includegraphics[width=0.33\textwidth]{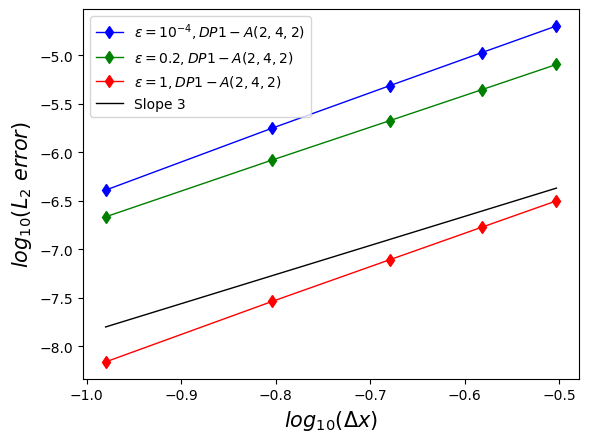}}
    \subfloat[CK-ARS WP]{\label{CKARSwpspace}\includegraphics[width=0.33\textwidth]{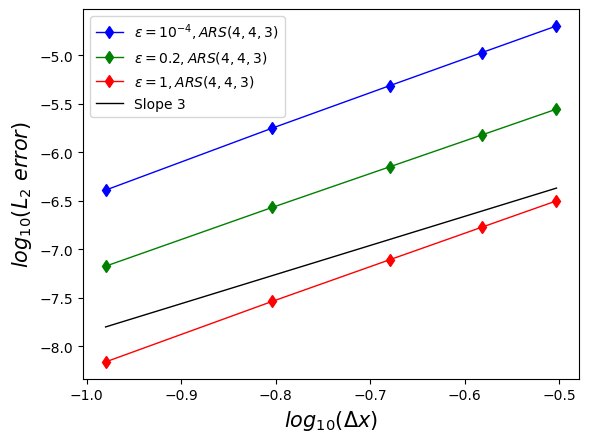}}
    \caption{Accuracy in space for the third order spatial scheme coupled with DP1-A$(2,4,2)$ (left) and ARS$(4,4,3)$ (right) for the time approximation. }
    \label{space order}
\end{figure}

The problem set-up is the same as described in the previous subsection, except for the following changes. Here, we consider the final time to be $T=0.01$ and $\Delta t=0.001$ 
so that the error in time is small enough to study the spatial accuracy. To do so, we consider the following number of points in space: $N_x=20, 24, 30, 40$ and $60$. The reference solution is obtained with $N_x=120$. \\
Since the spatial accuracy plots obtained from different time integrators are quite similar, we present only the plots obtained by using DP1-A$(2,4,2)$ and ARS$(4,4,3)$ for different values of $\epsilon$ $(\epsilon=10^{-4},0.2,1)$ in \cref{Anwpspace,CKARSwpspace}. For the spatial discretization, we only show the results obtained by 
the third order spatial matrices on non-staggered grid presented in \cref{Subsec: Non-stag grid} so that the scheme is expected to be third order accurate in space. 
In \cref{Anwpspace,CKARSwpspace}, the expected order is observed for the two time integrators and for the three considered values of $\epsilon$. 

\subsubsection{Qualitative results}
\begin{figure}
    \centering
    \subfloat[$\epsilon=1$]{\label{QualAem0}\includegraphics[width=0.33\textwidth]{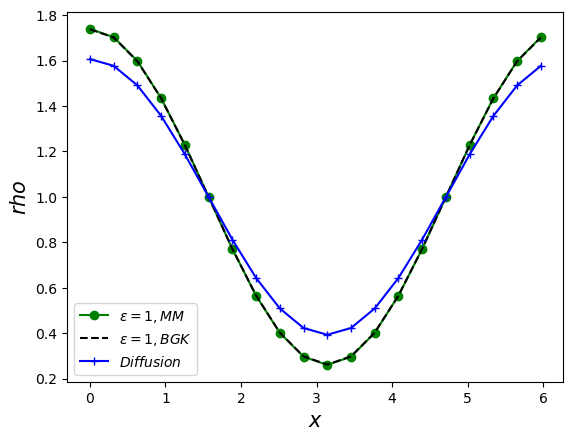}}
    \subfloat[$\epsilon=0.4$]{\label{QualAem1}\includegraphics[width=0.33\textwidth]{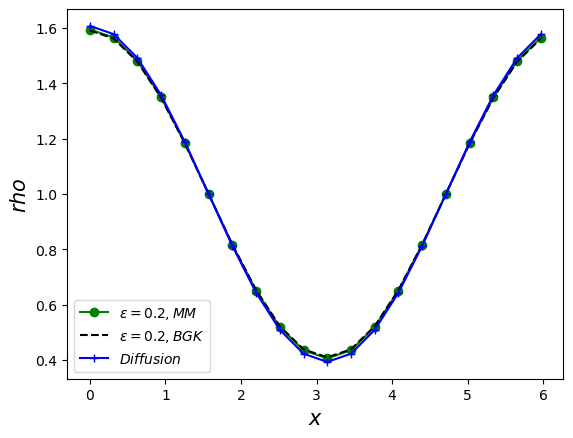}}
    \subfloat[$\epsilon=10^{-4}$]{\label{QualAem4}\includegraphics[width=0.33\textwidth]{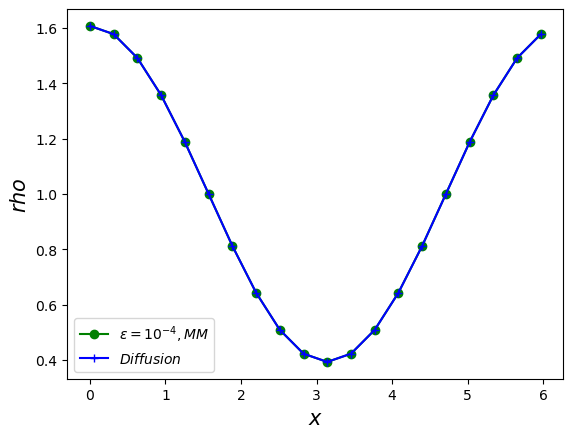}}
    \caption{Qualitative results for diffusion asymptotics}
    \label{Qual results}
\end{figure}

In this part, we compare the density obtained by the micro-macro equation (MM), the linear kinetic equation with BGK collision operator (BGK) and the asymptotic diffusion equation, for different values of $\epsilon$. 
The MM scheme described in previous sections is utilised,  the BGK is discretized using an IMEX (implicit treatment of collision term and explicit treatment of transport term) scheme whereas for the diffusion model, an implicit scheme is used. For all three models, the Butcher tableau corresponding to DP1-A$(2,4,2)$ time integrator is used. For the spatial discretization, we use third order scheme on non-staggered grid.\\
The problem domain $L=[0,2\pi]$ is discretised using $N_x=20$ grid points for all the three models. The final time is $T= 0.5$, and $\Delta t= 0.005$. We use the same N-WP initial and boundary conditions described in the previous 
subsection. Further, we also consider the same velocity discretization as before for both MM and BGK models.\\
In \cref{QualAem0} for rarefied regime ($\epsilon=1$), the MM and BGK models compare very well, while the diffusion model is different as expected. In the intermediate regime ($\epsilon= 0.2$), the BGK and MM models match very well while the diffusion model is slightly different. 
For $\epsilon=10^{-4}$, we only compare MM and the diffusion in \cref{QualAem4} and illustrate the AP property of the time integrators used for MM.


\subsection{Advection-diffusion asymptotics}
\begin{figure}
    \centering
    \subfloat[A]{\label{AdvDiffA}\includegraphics[width=0.33\textwidth]{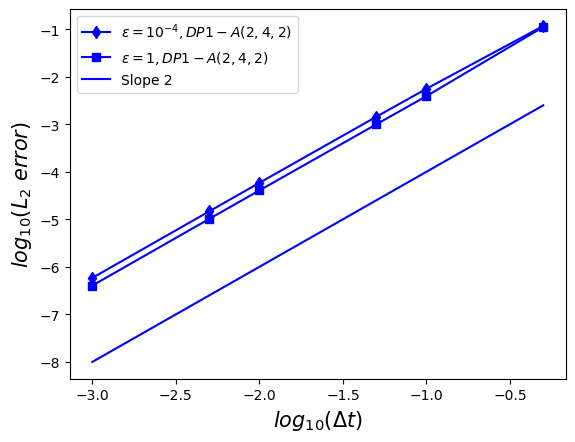}}
    \subfloat[CK-ARS]{\label{AdvDiffCKARS}\includegraphics[width=0.33\textwidth]{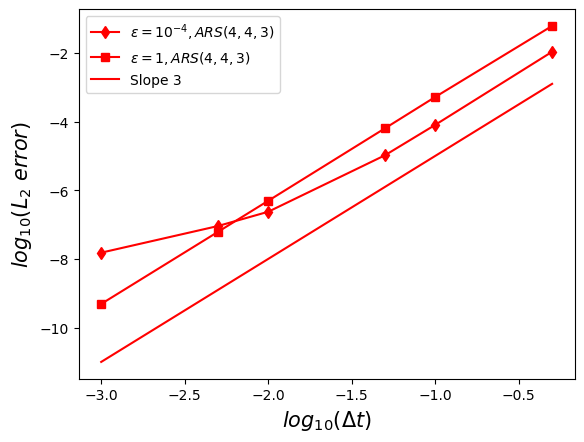}}
    \caption{Accuracy in time. Left:  DP1-A$(2,4,2)$  (N-WP initial data). Right: ARS$(4,4,3)$ (WP initial data). The reference solution is obtained from the micro-macro scheme with $\Delta t=10^{-4}$. }
    \label{time order AdvDiff}
\end{figure}

In this subsection, we present the time accuracy of our high order micro-macro scheme for the advection-diffusion case.
As in the diffusion case, the spatial domain $L=[0,2\pi]$ is discretised using $N_x=20$ grid points whereas the  velocity domain is $[-v_{\max}, v_{\max}]$ with $v_{\max}=5$ and $\Delta v=1$. The initial condition for the problem is: 
\begin{gather}
    \rho(0,x)=\sin(x) \\
    \text{Well-prepared data (WP): } g(0,x,v) = \epsilon^2 (I-\Pi) \left( v^2M \right) \rho (0,x) \\
    \text{Non-well prepared data (N-WP): } g(0,x,v) = (I-\Pi) \left( v^2M \right) \rho (0,x), 
\end{gather}
with $M(v)=\frac{1}{\sqrt{2\pi}} e^{-v^2/2}$. 
Periodic boundary conditions are used on both $\rho$ and $g$. The spatial terms are discretised by using the atmost-first order accurate matrices on staggered grid presented in \cref{Subsec: Stag grid}. 
The final time is $T=0.5$, and the following time steps  are considered: $\Delta t\!=\!0.5, 0.1, 0.05, 0.01, 0.005, 0.001$. We observe the time order of accuracy for both $\epsilon=1$ and $\epsilon=10^{-4}$. We choose the highest order time integrator in both type A and CK-ARS schemes for studying the time accuracy. Hence, we consider DP1-A$(2,4,2)$ and ARS$(4,4,3)$ with N-WP and WP data respectively. \\
Asymptotically, our micro-macro scheme degenerates to a consistent scheme for the advection-diffusion equation with advection and diffusion terms being treated explicitly and implicitly respectively. Hence, unlike the case of diffusion asymptotics for which an extra order is observed asymptotically, DP1-A$(2,4,2)$ remains second order accurate for $\epsilon=10^{-4}$ since both explicit and implicit matrices of the Butcher tableau are involved here (\cref{AdvDiffA}). For $\epsilon=1$, the required second order accuracy is observed. Further, the required third order accuracy of ARS$(4,4,3)$ is observed for both $\epsilon=10^{-4},1$ in \cref{AdvDiffCKARS}, since well-prepared initial data is considered. 

\subsection{Inflow boundary condition}
In this subsection, the high order numerical scheme for micro-macro model that allows inflow boundary conditions is validated numerically. We first present the time accuracy results for  high order schemes. Then, some  qualitative plots are shown for two tests with zero inflow at the right boundary, and equilibrium and non-equilibrium inflows respectively at the left boundary.   

\subsubsection{Time order of accuracy}
\begin{figure}
    \centering
    \subfloat[Type A, $\epsilon=1$ ]{\label{BdyAem0}\includegraphics[width=0.33\textwidth]{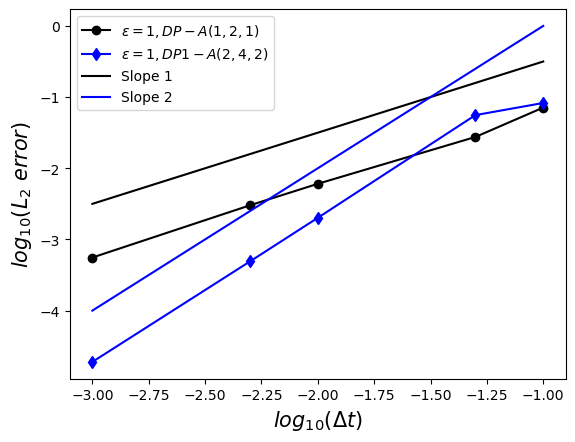}}
    \subfloat[Type A, $\epsilon=10^{-4}$]{\label{BdyAem4}\includegraphics[width=0.33\textwidth]{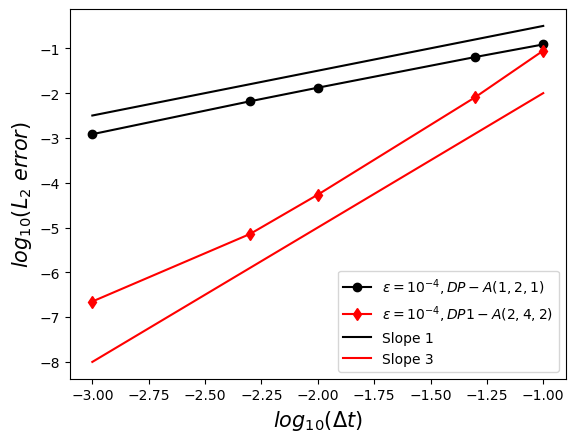}}

    
    \caption{Accuracy in time with type A schemes for $\epsilon=1$ (left) and $\epsilon=10^{-4}$ (right). The reference solution is obtained from the micro-macro for inflow boundaries scheme with $\Delta t=10^{-4}$. }
    \label{time order inflow}
\end{figure}

If the domain of the problem is a half-plane, $\omega(x,v)=\begin{bmatrix}-v,0,0,\cdots \end{bmatrix}$ can be chosen $\forall x$ as described in \cite{doi:10.1137/120865513}. Here, for numerical purposes, we consider a domain of $L=[0,2]$ and assume that the right boundary does not influence the dynamics. \\
The spatial domain is discretised using $N_x=20$ grid points and the velocity domain is $[-v_{\max}, v_{\max}]$ 
with $v_{\max}=5$ with $\Delta v=1$. 
The initial conditions at all interior points and right boundary conditions for the variables $\rho, \overline{\rho}$ and $\overline{g}$ are considered to be $0$. The left boundary conditions (for $v_k>0$) are:
\begin{equation}
    f\left(t,x_i=0,v_k\right) = M(v_k), \;\;\;     \overline{\rho} \left(t,x_i=0\right) = 1, \quad \overline{g}(t,x_{i+1/2}=-\Delta x/2,v_k) = 0,  
\end{equation}
with $M(v)=\frac{1}{\sqrt{2\pi}} e^{-v^2/2}$. 
The final time is $T=0.1$, and the following time steps are considered to check the accuracy in time: $\Delta t=0.1, 0.05, 0.01, 0.005, 0.001$. Like in the previous problems, we observe the time order of accuracy for both $\epsilon=1$ and $\epsilon=10^{-4}$. The time integrators considered are DP-A$(1,2,1)$ and DP1-A$(2,4,2)$. The reference solution for each curve in \cref{time order inflow} is obtained by using the same micro-macro scheme corresponding to that curve with  $\Delta t = 10^{-4}$. For type A time integrators with $\epsilon=1$ in \cref{BdyAem0}, first and second order accuracies of DP-A$(1,2,1)$ and DP1-A$(2,4,2)$ are observed. In \cref{BdyAem4} for $\epsilon=10^{-4}$, first and third order accuracies of DP-A$(1,2,1)$ and DP1-A$(2,4,2)$ respectively are observed. As for the (periodic) diffusion case, DP1-A$(2,4,2)$ turns out to be third order accurate since only the implicit part of Butcher tableau is involved asymptotically. 
For ARS$(2,2,2)$ and ARS$(4,4,3)$ time integrators (not shown here), 
order reduction to first order for $\epsilon=1$ (due to the initial condition). However, for $\epsilon=10^{-4}$, the required second and third orders respectively are observed. 

\subsubsection{Qualitative results for equilibrium inflow}
\begin{figure}
    \centering
    \subfloat[$\epsilon=1$]{\label{BdyEqbmAem0}\includegraphics[width=0.33\textwidth]{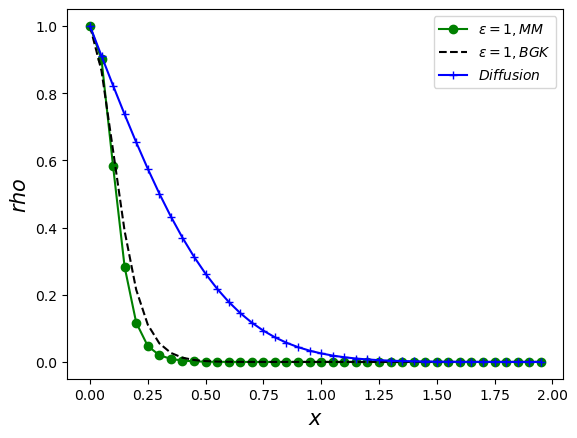}}
    \subfloat[$\epsilon=0.4$]{\label{BdyEqbmAem1}\includegraphics[width=0.33\textwidth]{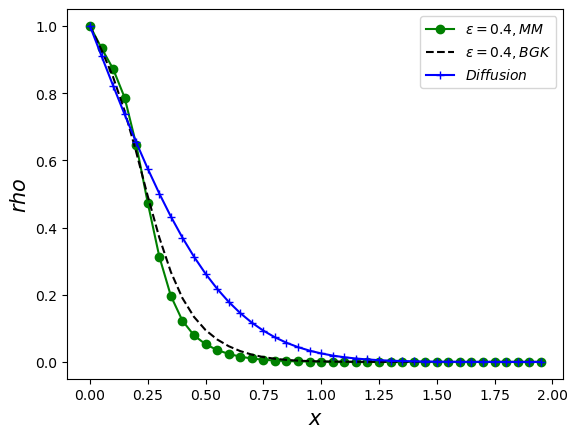}}
    \subfloat[$\epsilon=10^{-4}$]{\label{BdyEqbmAem4}\includegraphics[width=0.33\textwidth]{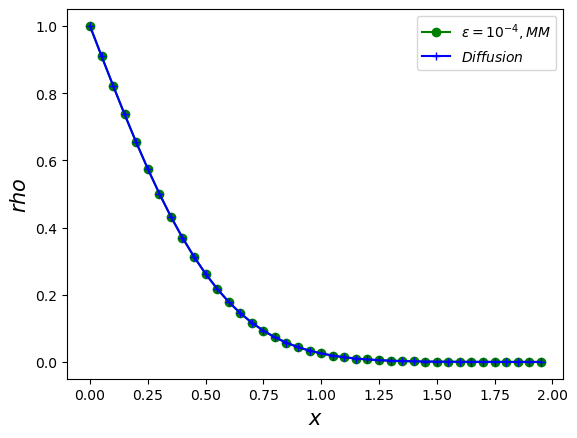}}
    \caption{Qualitative results for equilibrium inflow at the left boundary.}
    \label{BdyEqbmA}
\end{figure}

In this part, we consider the same problem as before and present a comparison of density plots obtained by using schemes based on micro-macro (MM), full-kinetic (BGK) and diffusion models, for different regimes of $\epsilon$. 
The  boundary conditions for diffusion model  $\rho(t, x=0)=1$ and $\rho(t, x=2)=0$. 
The final time is $T=0.1$, $N_x=40$ and $\Delta t =0.001$. Further, we consider the same velocity discretization as before for both MM and BGK models. The results for MM are obtained by DP1-A$(2,4,2)$ time integrator. \\
In \cref{BdyEqbmAem0} for rarefied regime ($\epsilon=1$), the MM and BGK results are in good agreement. In the intermediate regime ($\epsilon=0.4$) in \cref{BdyEqbmAem1}, the MM and BGK results are still close, and still  different from the diffusion one. For $\epsilon=10^{-4}$, only MM and the diffusion are plotted and are found to be in very good agreement, thereby illustrating the AP property of the numerical scheme for MM. 

\subsubsection{Qualitative results for non-equilibrium inflow}
\begin{figure}
    \centering
    \subfloat[$\epsilon=1$]{\label{BdyNonEqbmAem0}\includegraphics[width=0.33\textwidth]{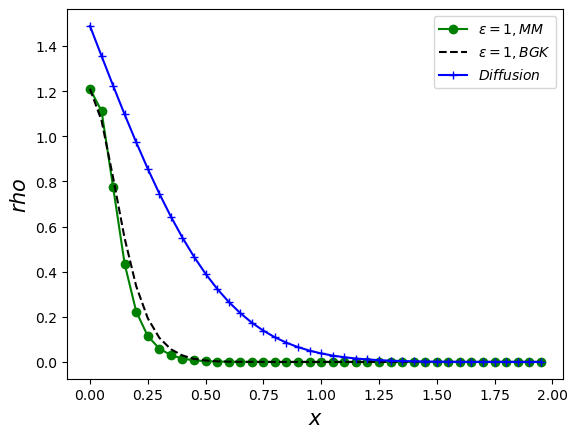}}
    \subfloat[$\epsilon=0.4$]{\label{BdyNonEqbmAem1}\includegraphics[width=0.33\textwidth]{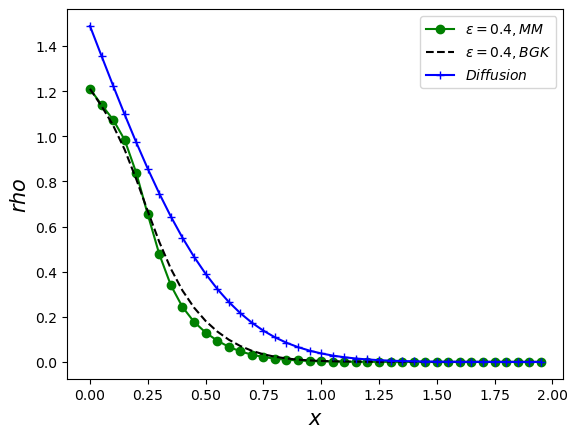}}
    \subfloat[$\epsilon=10^{-4}$]{\label{BdyNonEqbmAem4}\includegraphics[width=0.33\textwidth]{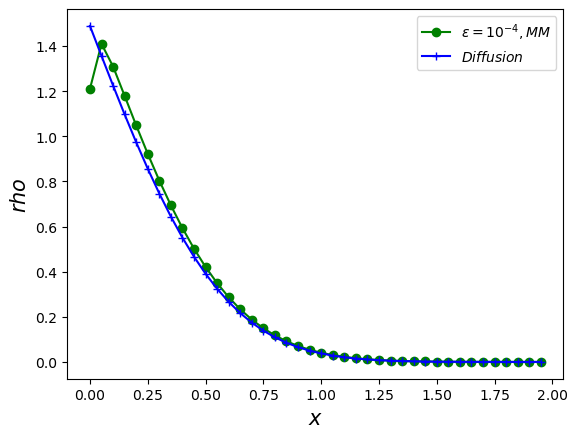}}
    \caption{Qualitative results  for non-equilibrium  inflow at the left boundary. }
    \label{BdyNonEqbmA}
\end{figure}

In this part, we consider the same problem as before, but the left boundary condition is chosen as (for $v_k>0$)
\begin{gather}
    f\left(t,x_i\!=\!0,v_k\right) = v_k M_k, \;\;\;\;  
    \overline{\rho} \left(t,x_i\!=\!0\right) = \left< f\left(t,x_i\!=\!0,v_k\right)\right>_{V_-} \\
    \overline{g}\bigl(t,x_{i+1/2}\!=\!-\tfrac{\Delta x}{2},v_k\bigr) \!=\! 2 \bigl( f\left(t,x_i\!=\!0,v_k\right) \!-\!  \overline{\rho} \left(t,x_i\!=\!0\right) M_k \bigr) \!-\! \overline{g}\bigl(t,x_{i+1/2}\!=\!\tfrac{\Delta x}{2},v_k\bigr). 
\end{gather}
The number of grid points, velocity discretization, final time and time step are the same as in the previous (equilibrium inflow) case. 
Here, we present a comparison of plots obtained by using schemes based on MM, BGK and diffusion models, for different regimes of $\epsilon$. The scheme described in subsection \ref{subsec:boundary} is used for the micro-macro model and a standard BGK approximation where only inflow boundary condition is needed serves as a reference. For diffusion, the diffusion term is treated implicitly and the left boundary condition for diffusion model is obtained from \cite{klar1} which translates in our context
\begin{multline*}
    \rho(t,x_i=0) = \frac{\sum_{v_k>0} v_k f\left(t,x_i=0,v_k\right) \Delta v}{\sum_{v_k>0} v_k M_k \Delta v} \\ + \frac{1}{ \kappa \sum_{v_k} M_k \Delta v} \sum_{v_k>0} v_k^2 \left( f\left(t,x_i=0,v_k\right) - M_k \frac{\sum_{v_k>0} v_k f\left(t,x_i=0,v_k\right) \Delta v}{\sum_{v_k>0} v_k M_k \Delta v} \right) \Delta v.
\end{multline*}
In \cref{BdyNonEqbmAem0} for rarefied regime ($\epsilon=1$), the MM and BGK models compare very well, while the diffusion model is driven by the macro boundary condition. In the intermediate regime ($\epsilon=0.4$) in \cref{BdyNonEqbmAem1}, in the MM and BGK results (which are in a good agreement), a boundary layer starts to be created whereas 
it is not the case for the diffusion model. For $\epsilon=10^{-4}$,  one can see that MM model develops a boundary layer at the left boundary before aligning with the diffusion model in the interior of the domain.  This is consistent with the results observed in the literature \cite{klar1, doi:10.1137/120865513, doi:10.1137/07069479X, crous_lemou}. 

\begin{appendix}
\section{Appendix: Matrix notation}
\label{Sec: App Matrix notation}
The $\mathsf{circ}$ function is given by:
\begin{gather}
    \mathsf{circ}([a_1,a_2,..,\underline{a_m},..,a_M])= 
    \left[\begin{smallmatrix} a_m & a_{m+1} & .. & a_M & 0 & .. & 0 & a_1 &  .. & a_{m-1} \\  a_{m-1} & a_m & a_{m+1} & .. & a_M & 0 & .. & 0 & a_1 &  .. \\ \\ \\ a_{m+2} & .. & a_M & 0 & .. & 0 & a_1 & .. & a_{m} & a_{m+1} \\ a_{m+1} & .. & a_M & 0 & .. & 0 & a_1 & .. & a_{m-1} & a_{m}  \end{smallmatrix}\right]
\end{gather}
The $\mathsf{circ_b}([-1,\underline{1}])_{(N_x-1)\times(N_x-2)}$ function is given by:
\begin{gather}
    \mathsf{circ_b}([-1,\underline{1}])_{(N_x-1)\times(N_x-2)}= \left[\begin{smallmatrix} 1 & 0 & .. & 0 \\  -1 & 1 & 0 & .. \\  \\ \\ .. & .. & -1 & 1 \\  .. & .. & .. & -1 \end{smallmatrix}\right]_{(N_x-1)\times(N_x-2)}
\end{gather}

\section{Appendix: Butcher tableau}
\label{App: Butcher tableau}
The following is the 2-stage second order accurate Butcher tableau ARS$(2,2,2)$:  
\begin{equation*}
\begin{tabular}{p{0.25cm}|p{0.4cm}p{0.7cm}p{0.5cm}}
\centering $0$ & $0$ & $0$ & $0$ \cr
\centering $\gamma$ & $\gamma$ & $0$ & $0$ \cr
\centering $1$ & $\delta$ & $1-\delta$ & $0$ \cr
\hline
 & $\delta$ & $1-\delta$ & $0$
\end{tabular}
 \quad \quad
\begin{tabular}{p{0.25cm}|p{0.4cm}p{0.7cm}p{0.5cm}}
\centering $0$ & $0$ & $0$ & $0$ \cr
\centering $\gamma$ & $0$ & $\gamma$ & $0$ \cr
\centering $1$ & $0$ & $1-\gamma$ & $\gamma$ \cr
\hline
 & $0$ & $1-\gamma$ & $\gamma$
\end{tabular}
\end{equation*}
Here, $\gamma=1-\frac{1}{\sqrt{2}}$ and $\delta=1-\frac{1}{2\gamma}$. \\
The following is the 4-stage third order accurate Butcher tableau ARS$(4,4,3)$:  
\begin{equation*}
\begin{tabular}{p{0.55cm}|p{0.9cm}p{0.8cm}p{0.55cm}p{0.8cm}p{0.2cm}}
\centering $0$ & $0$ & $0$ & $0$ & $0$ & $0$ \cr
\centering $1/2$ & $1/2$ & $0$ & $0$ & $0$ & $0$ \cr
\centering $2/3$ & $11/18$ & $1/18$ & $0$ & $0$ & $0$ \cr
\centering $1/2$ & $5/6$ & $-5/6$ & $1/2$ & $0$ & $0$ \cr
\centering $1$ & $1/4$ & $7/4$ & $3/4$ & $-7/4$ & $0$ \cr
\hline
 & $1/4$ & $7/4$ & $3/4$ & $-7/4$ & $0$
\end{tabular}
 \quad \quad
\begin{tabular}{p{0.55cm}|p{0.2cm}p{0.8cm}p{0.8cm}p{0.55cm}p{0.55cm}}
\centering $0$ & $0$ & $0$ & $0$ & $0$ & $0$ \cr
\centering $1/2$ & $0$ & $1/2$ & $0$ & $0$ & $0$ \cr
\centering $2/3$ & $0$ & $1/6$ & $1/2$ & $0$ & $0$ \cr
\centering $1/2$ & $0$ & $-1/2$ & $1/2$ & $1/2$ & $0$ \cr
\centering $1$ & $0$ & $3/2$ & $-3/2$ & $1/2$ & $1/2$  \cr
\hline
 & $0$ & $3/2$ & $-3/2$ & $1/2$ & $1/2$
\end{tabular}
\end{equation*}
 For type A, we use 2-stage first order accurate Butcher tableau DP-A$(1,2,1)$ ($\gamma \geq \frac{1}{2}$)\\
 \begin{equation*}
\begin{tabular}{p{0.25cm}|p{0.4cm}p{0.4cm}}
\centering $0$ & $0$ & $0$ \cr
\centering $1$ & $1$ & $0$  \cr
\hline
& $1$ & $0$
\end{tabular}
 \quad \quad
\begin{tabular}{p{0.25cm}|p{0.7cm}p{0.4cm}}
\centering $\gamma$ & $\gamma$ & $0$ \cr
\centering $1$ &  $1-\gamma$ & $\gamma$ \cr
\hline
 &  $1-\gamma$ & $\gamma$
\end{tabular}
\end{equation*}
The following is the 4-stage second order accurate Butcher tableau DP2-A$(2,4,2)$:
\begin{equation*}
\begin{tabular}{p{0.4cm}|p{0.4cm}p{0.7cm}p{0.7cm}p{0.4cm}}
\centering $0$ & $0$ & $0$ & $0$ & $0$  \cr
\centering $0$ & $0$ & $0$ & $0$ & $0$  \cr
\centering $1$ & $0$ & $1$ & $0$ & $0$  \cr
\centering $1$ & $0$ & $1/2$ & $1/2$ & $0$  \cr
\hline
 & $0$ & $1/2$ & $1/2$ & $0$
\end{tabular}
 \quad \quad
\begin{tabular}{p{0.4cm}|p{0.5cm}p{1.1cm}p{1.1cm}p{0.4cm}}
\centering $\gamma$ & $\gamma$ & $0$ & $0$ & $0$  \cr
\centering $0$ & $-\gamma$ & $\gamma$ & $0$ & $0$  \cr
\centering $1$ & $0$ & $1-\gamma$ & $\gamma$ & $0$  \cr
\centering $1$ & $0$ & $1/2$ & $1/2-\gamma$ & $\gamma$  \cr
\hline
 & $0$ & $1/2$ & $1/2-\gamma$ & $\gamma$
\end{tabular}
\end{equation*}
The following is the 4-stage second order accurate Butcher tableau DP1-A$(2,4,2)$ which achieves third order accuracy on the DIRK part:
\begin{equation*}
\begin{tabular}{p{0.6cm}|p{0.7cm}p{0.5cm}p{0.7cm}p{0.5cm}}
\centering $0$ & $0$ & $0$ & $0$ & $0$  \cr
\centering $1/3$ & $1/3$ & $0$ & $0$ & $0$  \cr
\centering $1$ & $1$ & $0$ & $0$ & $0$  \cr
\centering $1$ & $1/2$ & $0$ & $1/2$ & $0$  \cr
\hline
 & $1/2$ & $0$ & $1/2$ & $0$
\end{tabular}
 \quad \quad
\begin{tabular}{p{0.6cm}|p{0.9cm}p{1.1cm}p{0.7cm}p{0.6cm}}
\centering $1/2$ & $1/2$ & $0$ & $0$ & $0$  \cr
\centering $2/3$ & $1/6$ & $1/2$ & $0$ & $0$  \cr
\centering $1/2$ & $-1/2$ & $1/2$ & $1/2$ & $0$  \cr
\centering $1$ & $3/2$ & $1-3/2$ & $1/2$ & $1/2$  \cr
\hline
 & $3/2$ & $1-3/2$ & $1/2$ & $1/2$
\end{tabular}
\end{equation*}

\end{appendix}

\section*{Acknowledgement}
Megala Anandan sincerely acknowledges the indispensable support provided by Prof. S. V. Raghurama Rao, whose encouragement was instrumental in the successful completion of this research work.

\bibliographystyle{siamplain}
\bibliography{references}  
\end{document}